\documentclass{amsart}
\usepackage{amssymb,amsmath}

\theoremstyle{plain}
\newtheorem{theorem}{Theorem}
\newtheorem{lemma}[theorem]{Lemma}
\newtheorem{conjecture}[theorem]{Conjecture}
\newtheorem{proposition}[theorem]{Proposition}

\theoremstyle{definition}
\newtheorem{definition}{Definition}

\theoremstyle{remark}

 \newcommand\Ad{{\rm Ad}}
 \newcommand\rk{{\rm rk}}
 \newcommand\tr{{\rm tr}}

\def\mySp{{\rm Sp}}

\newcommand\Ca{\mathbb{O}}
\newcommand\C{\mathbb{C}}
\newcommand\LA{{\rm A}}

\newcommand\LC{{\rm C}}

\newcommand\LE{{\rm E}}
\newcommand\LF{{\rm F}}
\newcommand\LG{{\rm G}}

\newcommand\Le{\mathfrak e}

\newcommand\Lg{\mathfrak g}
\newcommand\N{{\rm N}}
\newcommand\R{\mathbb{R}}
\newcommand\SO{{\rm SO}}
\newcommand\SUxU[2]{ {\rm S(U(}#1) {\cdot} {\rm U(}#2{\rm ))} } 
\newcommand\SU{{\rm SU}}

\newcommand\Spin{{\rm Spin}}
\newcommand\T{{\rm T}}
\newcommand\U{{\rm U}}

\newcommand\act{{\,\mbox{\raisebox{0.2em}{\bf .}}\,}}

\newcommand\eK{{\e K}}
\newcommand\e{{\rm e}}

\newcommand\g{{\mathfrak g}}
\newcommand\hl{\rule[-.48em]{0em}{1.5em}}
\newcommand\h{{\mathfrak h}}
\newcommand\id{{\rm id}}
\newcommand\m{{\mathfrak m}}
\newcommand\pp{{\mathfrak p}}

\newcommand\p{{\mathfrak p}}

\newcommand\spin{\mathfrak s \mathfrak p \mathfrak i \mathfrak n}

\newcommand\str{\rule[-.48em]{0em}{1.7em}}

\newcommand\su{\mathfrak s \mathfrak u}
\newcommand\s{\sigma}

\newcommand{\lp}{\left (}                     
\newcommand{\rp}{\right )}                    
\newcommand{\lbc}{\left \{}                   
\newcommand{\rbc}{\right \}}                  

\begin{document}


\newcounter{refcounter}
\renewcommand\therefcounter{{\bf (\arabic{refcounter})}}
\newcommand\mylabel[1]{\refstepcounter{refcounter}\bf\therefcounter\label{#1}}

%
%
%

\title[Low cohomogeneity and polar actions]
{Low cohomogeneity and polar actions on exceptional compact Lie groups}

\author{Andreas Kollross}

\address{Institut f\"{u}r Mathematik\\Universit\"{a}t Augsburg\\86135
Augsburg\\Germany}

\email{kollross@math.uni-augsburg.de}


\maketitle

\begin{abstract}
We study isometric Lie group actions on the compact exceptional groups $\LE_6$,
$\LE_7$, $\LE_8$, $\LF_4$ and $\LG_2$ endowed with a biinvariant metric. We
classify polar actions on these groups, in particular, we show that all polar
actions are hyperpolar. We determine all isometric actions of cohomogeneity
less than three on $\LE_6$, $\LE_7$, $\LF_4$ and all isometric actions of
cohomogeneity less than~$20$ on~$\LE_8$. Moreover we determine the principal
isotropy algebras for all isometric actions on~$\LG_2$.
\end{abstract}


\section*{Introduction}
\renewcommand\thesubsection{\arabic{subsection}}

\label{Intro}


We study isometric Lie group actions on the compact exceptional simple Lie
groups $\LE_6$, $\LE_7$, $\LE_8$, $\LF_4$ and $\LG_2$ endowed with a
biinvariant Riemannian metric; we classify actions with low cohomogeneity and
polar actions on these spaces. An isometric action of a compact Lie group on a
Riemannian manifold is called {\em polar} if there exists an immersed connected
submanifold~$\Sigma$ which intersects the orbits orthogonally and meets every
orbit. Such a submanifold~$\Sigma$ is called a {\em section} of the group
action. If the section is flat in the induced metric, the action is called {\em
hyperpolar}. Our main result is the following.

\begin{theorem}\label{ThmPolar}
Let $L$ be a connected simple compact Lie group of type $\LE_6$,
$\LE_7$, $\LE_8$, $\LF_4$ or $\LG_2$ endowed with a biinvariant Riemannian
metric. Let $H \subseteq L \times L$ be a closed subgroup such that the action
of $H$ on~$L$ defined as in~(\ref{BiAction}) is polar. Then the $H$-action
on~$L$ is hyperpolar or the $H$-orbits are finite.

Moreover, if the cohomogeneity of the polar $H$-action on~$L$ is greater than
two, then $H$ is a symmetric subgroup of $L \times L$.
\end{theorem}

In the course of proving Theorem~\ref{ThmPolar}, we obtain an explicit
description of all polar actions of connected groups on the exceptional compact
Lie groups. As a further result, we classify actions of certain low
cohomogeneities on the exceptional groups, cf.\ Theorems~\ref{ThmCohOne},
\ref{ThmCohTwo}, \ref{ThmPrOrbG2}, \ref{ThE8Low}.


It should be noted that the classification problem for polar actions in the
special case that the section is flat, i.e.\ for hyperpolar actions, had
been solved before. In fact, the author has classified hyperpolar actions on
all irreducible compact symmetric spaces in~\cite{hyperpolar}.

If the additional assumption that the section is flat is dropped, i.e.\ if one
considers actions on irreducible compact symmetric spaces which are polar, but
not necessarily hyperpolar, then there is a sharp contrast between the case of
rank-one symmetric spaces and the higher rank case; while there are many
examples of polar actions with non-flat sections on rank-one symmetric spaces,
see~\cite{pth1} for a classification, there are as yet no examples known on the
spaces of higher rank. In fact, there is the following conjecture.

\begin{conjecture}[Biliotti~\cite{biliotti}]
Any polar action with orbits of positive dimension on an irreducible compact
Riemannian symmetric space of higher rank is hyperpolar.
\end{conjecture}

This conjecture was shown to be true for all symmetric spaces of type~I, i.e.\
symmetric spaces $G/K$ where $G$ is a simple compact Lie group and $K$ is a
symmetric subgroup, by the author in~\cite{polar}.  Earlier the conjecture had
been proved to hold for actions with a fixed point by Br\"{u}ck~\cite{brueck},
for actions on the complex quadrics by Podest\`{a} and Thorbergsson~\cite{pt}, on
complex Grassmannians by Biliotti and Gori~\cite{bg}, and by
Biliotti~\cite{biliotti} for compact irreducible Hermitian symmetric spaces. It
remained open for the case of symmetric spaces of type~II, i.e.\ the simple
compact connected Lie groups equipped with a biinvariant metric.

Our Theorem~\ref{ThmPolar} now confirms Biliotti's conjecture in the special case
of exceptional compact Lie groups. However, the conjecture still remains open for polar
actions on the classical Lie groups $\SO(n)$, $\SU(n)$ and $\mySp(n)$.


A prominent special case of (hyper)polar actions of independent interest, which has been studied by
many authors, is the case of {\em cohomogeneity one actions}, i.e.\ such actions
where the principal orbits are hypersurfaces. Hsiang and Lawson~\cite{hsl},
Takagi~\cite{takagi}, D'Atri~\cite{datri}, and Iwata~\cite{iwata} have determined all
cohomogeneity one actions on~${\rm S}^n$, $\C {\rm P}^n$, ${\mathbb{H}} {\rm P}^n$ and $\Ca {\rm P}^2$, respectively.

In \cite{hyperpolar} the author has classified cohomogeneity one actions on
all irreducible compact symmetric spaces, in particular, on simple compact Lie
groups. However, the classification there is only up to orbit equivalence.
Theorem~\ref{ThmCohOne} is therefore a refinement of this classification in
that all connected closed subgroup of the isometry group are given which
act with cohomogeneity one.

Motivated by the interest in cohomogeneity one actions, we carry on the study of actions whose principal orbits have low codimension in this article and classify actions of cohomogeneity two on the exceptional compact Lie groups, cf.\ Theorem~\ref{ThmCohTwo}.

For the groups $\LG_2$ and $\LE_8$, we can further improve these results. It turns out that with few exceptions, given by Theorem~\ref{ThmPrOrbG2}, all isometric actions on $\LG_2$ have finite principal isotropy groups and hence we have, in particular, determined the (co)dimensions of the principal orbits of all isometric actions on~$\LG_2$. Finally, we classify all isometric actions on~$\LE_8$ of cohomogeneity less than~$20$, see~Theorem~\ref{ThE8Low}.


This article is organized as follows. Theorem~\ref{ThmPolar} is proved in
Sections~\ref{Subactions}--\ref{Regular}. The rest of the article is concerned
with actions of low cohomogeneity. In Section~\ref{G2Isotropy} we determine the
Lie algebra type of the principal isotropy subgroups for every isometric action
of a compact Lie group on~$\LG_2$.
In Section~\ref{CohOneTwo}, we determine all isometric actions of compact Lie
groups on the exceptional groups where the cohomogeneity is less than three.
Since these actions occur as candidates for polar actions, we can use the proof
of Theorem~\ref{ThmPolar} to a large extent. In Section~\ref{LowCoh} we
classify low cohomogeneity actions on~$\LE_8$.


\section{Preliminaries}
\label{Prelim}


In this article, our objects of study are simple compact connected Lie
groups~$L$, endowed with a biinvariant Riemannian metric. Such a metric is
unique up to a constant scaling factor, whose choice is of course irrelevant
here; we may for instance assume that $L$ is equipped with the homogeneous
metric induced by the negative of the Killing form.

Now let $H$ be a compact Lie group acting isometrically on~$L$. The action is
polar if and only if the action restricted to the connected component of~$H$ is
polar and furthermore the cohomogeneity of the $H$-action remains the same if
the action is restricted to the connected component of~$H$. Therefore we will
assume that $H$ is a closed connected subgroup of~$L \times L$ (which covers
the connected component of the isometry group of~$L$) and that the action
of~$H$ on~$L$ is given by
\begin{equation}\label{BiAction}
(h_1,\,h_2) \act \ell = h_1 \, \ell \, h_2^{-1}\, \mbox{ for }\, (h_1,h_2) \in
H,\; \ell \in L.
\end{equation}

Assume the groups $H_1$ and $H_2$ act isometrically on the Riemannian manifolds
$M_1$ and $M_2$, respectively. The $H_1$-action on $M_1$ and the $H_2$-action
on $M_2$ are called {\em conjugate} if there exists an isometry $F \colon M_1
\to M_2$ and an isomorphism $\phi \colon H_1 \to H_2$ such that $F \lp g \act p
\rp = \phi(g) \act F(p)$ for all $g \in H$, $p \in M_1$. For the purposes of
this article, it obviously suffices to consider actions up to conjugacy.

Let $L$ be a semisimple compact Lie group equipped with the biinvariant metric
induced by the negative of the Killing form and let  $H$ be a closed subgroup
of $L \times L$. Then any automorphism $\sigma \colon L \to L$ is an isometry
and the $H$-action on~$L$ is conjugate to the action of~$\phi( H )$ on~$L$
where $F = \sigma$ and $\phi(h_1,h_2) = \lp \sigma(h_1) , \sigma(h_2) \rp$. Let
$\ell, r \in L$, then the map $L \to L$, $F \colon g \mapsto \ell \, g \,
r^{-1}$ is an isometry of~$L$ and the $H$-action is conjugate to the action of
$\lbc \lp \ell \, h_1 \, \ell^{-1}, r \, h_2 \, r^{-1} \rp \mid (h_1,h_2) \in H
\rbc$. Furthermore, the map $F \colon L \to L,\, g \mapsto g^{-1}$ is an
isometry and the $H$-action on~$L$ is conjugate to the action of $\lbc
(h_2,h_1) \in L \times L \mid (h_1,h_2) \in H \rbc$. However, notice that if
$\alpha$ is an outer automorphism of~$L$, then the action of $\lbc
(h_1,\alpha(h_2)) \in L \times L \mid (h_1,h_2) \in H \rbc$ is in general not
conjugate to the $H$-action on~$L$, see \cite{hyperpolar}, Theorem~3.2 and the
preceding remarks.

\paragraph{\bf Notation.}
Since in this article we frequently have to deal with certain subgroups of $L
\times L$, where $L$ is a simple compact Lie group, it is convenient to adopt
the following notational convention: If $H_1$ and $H_2$ are subgroups of~$L$,
then $H_1 \times H_2$ {\em always} denotes the subgroup $\lbc (h_1,h_2) \mid
h_1 \in H_1,\, h_2 \in H_2 \rbc$ of $L \times L$ and we avoid the use the
symbol ``$\times$'' altogether whenever we consider direct products of groups
otherwise. If $K$ is a subgroup of~$L$, we denote by $\Delta K$ the diagonally
embedded subgroup $\lbc (g,g) \mid g \in K \rbc$ of~$L \times L$. More
generally, if $K$ is a subgroup of~$L$ and $\sigma \colon L \to L$ is an
automorphism of~$L$, then $\Delta^{\sigma} K$ stands for the subgroup $\lbc
(g,\sigma(g)) \mid g \in K \rbc$ of~$L \times L$. Whenever we consider a closed
subgroup~$H$ of $L \times L$ given by one of the three notations just described
above, it is henceforth always understood that the action of~$H$ on~$L$ is
given by~(\ref{BiAction}). In some cases we use the notation of~\cite{dynkin1}
to uniquely describe the conjugacy class of a subgroup in an exceptional
compact Lie group, e.g.\ in Table~\ref{TMaxE6}, we denote by $\LG_2^1$ and
$\LG_2^3$ two non-conjugate subgroups of~$\LE_6$ which have Dynkin index $1$
and $3$, respectively.

If $H_1 \subset L$ is a closed subgroup such that its Lie algebra is the fixed
point set of an involutive automorphism of the Lie algebra of~$L$, then we call
$H_1$ a {\em symmetric subgroup} of~$L$. If $H = H_1 \times H_2$ where $H_1,
H_2 \subset L$ are symmetric subgroups, then the $H$-action on $L$ is called a
{\em Hermann action}~\cite{hermann}. The action of $\Delta^{\sigma} L$ on~$L$
is called the {\em $\sigma$-action} on~$L$. Hermann actions and
$\sigma$-actions are well known to be hyperpolar~\cite{hptt}.

If $H$ is a closed connected subgroup of~$L \times L$ and $H'$ is a closed
subgroup of~$H$ then we will refer to the action of~$H'$ on~$L$ as a {\em
subaction} of the $H$-action. If in addition the $H'$-action and the $H$-action
on~$L$ are orbit equivalent, we say the $H'$-action is an {\em orbit equivalent
subaction} of the $H$-action.


Let $L$ be a simple compact connected Lie group and let ${\mathfrak l}$ be its Lie
algebra. Let $H$ be a closed connected subgroup of~$L \times L$. Let $\pi_1$,
$\pi_2 \colon {\mathfrak l} \times {\mathfrak l} \to {\mathfrak l}$ be the canonical projections such that
$\pi_1(X,Y)=X$ and $\pi_2(X,Y)=Y$. Let $\h_1 = \ker \pi_2|_{\h}$, $\h_2 = \ker
\pi_1|_{\h}$. Since $\h$ is reductive, there exists a complementary ideal
$\h_3$ of $\h_1 + \h_2$ in~$\h$. Since $\h_1 \cap \h_2 = 0$, the Lie
algebra~$\h$ of~$H$ is isomorphic to the direct sum $\h_1 \oplus \h_2 \oplus
\h_3$ of Lie algebras.


\begin{lemma}\label{LmMaxSubgr}
Let $L$ be a simple compact connected Lie group and let $H$ be a closed
connected subgroup of~$L \times L$ such that the $H$-action on~$L$ is not
transitive. Then $H$ is contained in at least one of the following subgroups
of~$L \times L$.
\begin{enumerate}

\item $\Delta^{\sigma} L$, where $\sigma \colon L \to L$ is an automorphism
    of~$L$.

\item $H' \times H''$, where $H', H'' \subset L$ are maximal connected
    subgroups of~$L$.

\end{enumerate}
\end{lemma}


\begin{proof}
With the notation as above, consider $\h = \h_1 \oplus \h_2 \oplus \h_3$.
Assume that $\h_1 \cong {\mathfrak l}$. Then it follows that $H$ contains the subgroup $L
\times \lbc 1 \rbc$, which acts transitively on~$L$, a contradiction. The same
argument shows that $\h_2 \not\cong {\mathfrak l}$. Assume $\h_3 \cong {\mathfrak l}$. Since
$\pi_1|_{\h_1 \oplus \h_3}$ and $\pi_2|_{\h_2 \oplus \h_3}$ are injective, it
follows that $\h_1 = \h_2 = \lbc 0 \rbc$ and that $\pi_1|_{\h_3}$ and
$\pi_2|_{\h_3}$ are Lie algebra isomorphisms $\h_3 \to {\mathfrak l}$. It follows that $\h
= \lbc (\pi_1 (X), \pi_2 (X) \mid X \in \h_3 \rbc$ and thus $H =
\Delta^{\sigma} L$, where $\lp \sigma_* \rp_{\e} = \pi_2|_{\h_3} \circ
\pi_1|_{\h_3}^{-1}$. Now assume $\h_3 \not\cong {\mathfrak l}$. It follows that $\h
\subseteq \pi_1(\h) \times \pi_2(\h)$, where $\pi_i(\h)$ are proper subalgebras
of~${\mathfrak l}$.
\qed\end{proof}


The maximal connected subgroups of compact Lie groups are classified in
\cite{dynkin1} and \cite{dynkin2}, cf.\ also~\cite{hyperpolar}, Theorems~2.1
and~2.2. There is the following criterion for polarity of an isometric action
on a symmetric space. For a proof see~\cite{gorodski} or~\cite{hyperpolar}.
Note that sections of polar actions on Riemannian manifolds are always totally
geodesic submanifolds.


\begin{proposition}\label{PropPolCrit}
Let $G$ be a connected compact Lie group, let $K \subset G$ be a symmetric
subgroup and let $\g = {\mathfrak k} \oplus \pp$ be the Cartan decomposition. Let $H
\subseteq G$ be a closed subgroup. Let $k$ be the cohomogeneity of the
$H$-action on $G$. Then the following are equivalent.
\begin{enumerate}

\item The $H$-action on $G / K$ is polar w.r.t\ some Riemannian metric
    induced by an $\Ad(G)$-invariant scalar product on~$\g$.

\item For any $g \in G$ such that $gK$ lies in a principal orbit of the
    $H$-action on~$G / K$ the subspace $\nu = g^{-1} N_{g K}(H \act gK)
    \subseteq \pp$ is a $k$-dimensional Lie triple system such that the Lie
    algebra~${\mathfrak s} = \nu \oplus [\nu, \nu]$ generated by $\nu$ is orthogonal
    to $\Ad(g^{-1})\h$.

\item The normal space $\N_{\eK}(H \act \eK) \subseteq \pp$ contains a
    $k$-dimensional Lie triple system $\nu$ such that the Lie algebra~${\mathfrak s}
    = \nu \oplus [\nu, \nu]$ generated by $\nu$ is orthogonal to $\h$.

\end{enumerate}
\end{proposition}


\paragraph{\bf Remark.} Hyperpolar actions are characterized by the additional
property that the Lie triple system~$\nu$ in Proposition~\ref{PropPolCrit} is
abelian (in which case the Lie algebra~${\mathfrak s}$ is equal to $\nu$). In case the
$H$-action on $G/K$ is polar, the Lie triple system~$\nu$ corresponds to the
tangent space of a section containing $\eK$. Note that cohomogeneity one
actions are hyperpolar.

Let us now apply the criterion from Proposition~\ref{PropPolCrit} to the case
of a compact Lie group~$L$ equipped with a biinvariant metric. To this end, we
present the symmetric space~$L$ homogeneously as $G / K = (L \times L) / \Delta
L$. The Lie algebra of $K = \Delta L$ is
\begin{equation}
{\mathfrak k} = \left \{  \left ( X,\, X \right ) \mid X \in {\mathfrak l} \right \},
\end{equation}
where ${\mathfrak l}$ denotes the Lie algebra of~$L$. The Cartan complement of ${\mathfrak k}$ in $\g
\cong {\mathfrak l} \oplus {\mathfrak l}$ is
\begin{equation}
\p = \left \{  \left ( X,\, -X \right ) \mid X \in {\mathfrak l} \right \}.
\end{equation}

Assume the subgroup $H \subset G$ is of the form $H = H_1 \times H_2$, where
$H_1$ and $H_2$ are closed subgroups of $L$, i.e.\
\begin{equation}\label{SplitGroup}
H = \left \{  \left ( h_1,\, h_2 \right ) \mid h_1 \in H_1,\,h_2 \in H_2 \right
\}.
\end{equation}
Let $\m_1$ and $\m_2$ be the orthogonal complements of $\h_1$ and $\h_2$,
respectively in the Lie algebra~${\mathfrak l}$ of~$L$. By conjugation of $H$ we may
assume without loss of generality that the identity element~$\e$ of~$G$ lies in
a principal orbit of the $H$-action on $G / K$. Then the subspace $\nu \subset
\g$ in Proposition~\ref{PropPolCrit}~(ii) is given by
\begin{equation}
\nu = \left \{  \left ( X,\, -X \right ) \mid X \in \m_1 \cap \m_2 \right \}
\end{equation}
and $[\nu, \, \nu]$ is spanned by the elements
$$
\lp \left [ X,\, Y \right ],\left [ X,\, Y \right ] \rp,\quad X,\,Y \in \m_1
\cap \m_2.
$$
If now, say, $H_2$ is a symmetric subgroup of~$L$, then we have from the Cartan
relations $[ \m_2,\, \m_2] \subseteq \h_2$. Thus if $[\nu, \, \nu] \neq 0$ then
it follows that $[\nu, \, \nu]$ is not orthogonal to~$\h$. We have shown:

\begin{lemma}\label{LmNoLift}
Let $L$ be a compact Lie group equipped with a biinvariant metric. Let $H = H_1
\times H_2 \subset L \times L$ be closed subgroup as in~(\ref{SplitGroup}) and
such that $H_2 \subset L$ is a symmetric subgroup. Then the action of $H$ on
$L$ is polar if and only if it is hyperpolar.
\end{lemma}


The following Theorem was shown in~\cite{polar}, Theorem~5.4


\begin{theorem}\label{ThProductOfSpheres}
If a compact connected Lie group acts polarly and non-trivially on an
irreducible compact symmetric space then every section is covered by a
Riemannian product of spaces which have constant curvature.
\end{theorem}


\begin{lemma}\label{LmDimSection}
Let $L$ be a simple compact Lie group and let $H \subset L \times L$ be closed
subgroup action polarly on~$L$. Then $\dim H \ge \dim L - 3 \cdot \rk L$. In
particular, a compact connected nontrivial Lie group acting polarly on one of
the exceptional groups $\LG_2$, $\LF_4$, $\LE_6$, $\LE_7$, $\LE_8$ is of
dimension greater or equal $8$, $40$, $60$, $112$, $224$, respectively.
\end{lemma}


\begin{proof}
This is an immediate consequence of~\cite{polar}, Lemma~3.3.
\qed\end{proof}


\begin{lemma}\label{LmDimBd}
Let $L$ be an exceptional simple compact Lie group $\LF_4$, $\LE_6$, $\LE_7$,
or $\LE_8$ and let
$$
K = \left \{ (k_1, k_2) \mid k_1 \in K_1,\; k_2 \in K_2 \right \} \subseteq L
\times L,
$$
where $K_1$ and $K_2$ are closed proper subgroups of~$L$. If a closed subgroup
of $K$ acts polarly on~$L$, then $\dim K_i \ge 12$, $15$, $34$, $90$,
respectively.
\end{lemma}

\begin{proof}
Let $H \subseteq K$ be a closed connected subgroup acting polarly on~$L$. Let
$\h = \h_1 \oplus \h_2 \oplus \h_3$ be the Lie algebra of~$H$ with the notation
as in Lemma~\ref{LmMaxSubgr}. Then $\pi_2(\h) \cong \h_2 \oplus \h_3$ is a
subalgebra in the Lie algebra of~$K_2$. Assume $\h_2$ corresponds to a
symmetric subgroup of~$L$. Then $\pi_2(\h_2)$ a maximal subalgebra of~${\mathfrak l}$ and
it follows that $\h_3 = 0$. Thus the $H$-action is hyperpolar by
Lemma~\ref{LmNoLift} in this case and it follows that $\dim K_1 \ge 12$, $20,$
$47,$ $104$, respectively, see \cite{hyperpolar}, Section~2.4.4.

Now assume $\h_2$ does not correspond to a symmetric subgroup of~$L$. We
determine non-symmetric connected subgroups of maximal dimension in the groups
$\LF_4$, $\LE_6$, $\LE_7$ and $\LE_8$. Such groups are $\Spin(8) \subset
\LF_4$, $\Spin(10) \subset \LE_6$, $\LE_6 \subset \LE_7$ and $\LE_7 \cdot \U(1)
\subset \LE_8$, see~\cite{dynkin1}, cf.\ also Tables~\ref{TMaxF4}, \ref{TMaxE6}
and \ref{TMaxE7}. Thus the maximal dimension of a proper closed non-symmetric
subgroup in $\LF_4$, $\LE_6$, $\LE_7$, or $\LE_8$ is $28$, $45$, $78$ and
$134$, respectively and the assertion of the Lemma now follows directly from
Lemma~\ref{LmDimSection}, since $\pi_1(\h) \cong \h_1 \oplus \h_3$ is a
subalgebra in the Lie algebra of~$K_1$.
\qed\end{proof}

Assume a Lie group~$G$ acts isometrically on a Riemannian manifold~$M$ and let
$p \in M$. Then the {\em isotropy subgroup} $G_p = \lbc g \in G \,|\, g \act p
= p \rbc$ acts on~$\T_p M$ such that the tangent space~$\T_p \lp G \act p \rp$
and the normal space~$\N_p \lp G \act p \rp$ to the $G$-orbit through~$p$ are
invariant subspaces. The action of~$G_p$ on~$\N_p \lp G \act p \rp$ is called
the {\em slice representation} of the $G$-action at~$p$. The slice
representation is trivial if and only if $p$ lies in a principal orbit. Slice
representations are an import tool for analyzing Lie group actions on manifolds
since they provide a local linearization of the group action. In fact, they
will be our main tool in this article. A slice representation has the same
cohomogeneity as the action of~$G$ on~$M$; if the action of~$G$ on~$M$ is
polar, then all slice representations are polar, too; however, the converse is
not true, see Section~\ref{CohOneTwo}.


\section{Polarity of subactions}
\label{Subactions}


For hyperpolar actions on symmetric spaces of the compact type one has the
maximality property that the action of a closed subgroup~$H$ of the isometry
group can only have hyperpolar subactions if the $H$-action is itself
hyperpolar (maybe transitive); this follows immediately from
Proposition~\ref{PropPolCrit}. In fact, in case the symmetric space is
irreducible, the following stronger statement holds: If there is an inclusion
relation between two closed subgroups of the isometry group which both act
hyperpolarly then the actions are orbit equivalent (see below) or one action is
transitive, see Corollary~D of~\cite{hl}. Therefore, it is sufficient to
consider only maximal nontransitive subgroups of the isometry group in order to
find all groups acting hyperpolarly on a given irreducible compact symmetric
space, cf.\ the classification in~\cite{hyperpolar}.

However, for polar actions such a maximality property does not hold in general;
for example, let $G_1, G_2$ be compact Lie groups acting orthogonally on
$\R^{n_1}$ and $\R^{n_2}$, respectively, such that the $G_1$-action is polar,
but the $G_2$-action is not. Then the action of the direct product of $G_1$ and
$G_2$ on ${\rm S}^{n_1 + n_2-1} \subset \R^{n_1} \oplus \R^{n_2}$ is not polar,
even though the action restricted to the subgroup~$G_1$ is. To overcome this
difficulty, we introduce the notion of polarity minimality, which means that an
action does not have any polar subactions in a nontrivial way.

We say that an action of a group~$G$ on a Riemannian manifold~$M$ is {\em orbit
equivalent} to the action of a group~$G'$ on a Riemannian manifold~$M'$ if
there is an isometry $F \colon M \to M'$ such that $F$ maps any connected
component of a $G$-orbit in~$M$ onto a connected component of a $G'$-orbit
in~$M'$.


\begin{definition}
Let $G$ be a compact Lie group acting isometrically on a Riemannian manifold.
We say the action of~$G$ on~$M$ is {\em polarity minimal} if any closed
connected subgroup of~$G$ whose action on~$M$ is nontrivial and not orbit
equivalent to the $G$-action is non-polar.
\end{definition}


Note that a polarity minimal action can be polar or non-polar. We cite the
following proposition from~\cite{polar}, which gives some sufficient conditions
for an orthogonal representation to be polarity minimal.


\begin{proposition}\label{PolMinCriteria}
Let $\rho \colon G \to {\rm O}(V)$ be a representation of the compact connected Lie
group~$G$. Then $\rho$ is polarity minimal if one of the following holds.
\begin{enumerate}

\item[(i)] The representation $\rho$ is irreducible of cohomogeneity~$\ge
    2$.

\item[(ii)] The representation space~$V$ is the direct sum of two
    equivalent $G$-modules.

\item[(iii)] The representation space $V$ contains a $G$-in\-var\-iant
    sub\-module~$W$ such that the $G$-re\-pre\-sen\-ta\-tion on~$W$ is
    almost effective, non-polar, and polarity minimal.

\end{enumerate}
\end{proposition}


We will use Proposition~\ref{PolMinCriteria} to show that in many cases an
action on an exceptional group has a polarity minimal slice representation.
Under various conditions, some of which are collected in the following
proposition, this is sufficient to show that the action under consideration is
polarity minimal itself. This will be the main tool of our classification.


\begin{lemma}\label{LmPolMinHered}
Let $G$ be compact Lie group and $K \subset G$ be symmetric subgroup such that
$M = G / K$ is an irreducible symmetric space and let $H \subset G$ be a closed
subgroup. The action of $H$ on $M$ is non-polar and polarity minimal if there
is a non-polar polarity minimal submodule~$V \subseteq \N_p(H \act p)$ of the
slice representation at~$p$ such that one of the following holds.
\begin{enumerate}

\item $M$ is Hermitian symmetric and $\dim(V) > \rk(H)$.

\item $\dim(V) > s(M)$, where $s(M)$ is the maximal dimension of a totally
    geodesic submanifold of~$M$ locally isometric to a product of spaces
    with constant curvature, cf.~\cite{polar}, Lemma~3.3.

\item $V \subseteq \pp = \T_p M$ (where $\g = {\mathfrak k} \oplus \pp$ as usual such
    that ${\mathfrak k}$ is the Lie algebra of~$K = G_p$) contains a Lie triple system
    corresponding to an irreducible symmetric space of nonconstant
    curvature, e.g.\ an irreducible symmetric space of higher rank.

\item The isotropy group $H \cap K$ acts almost effectively on $V$ and
    $\rk(H \cap K) = \rk(H)$.

\end{enumerate}
\end{lemma}


\begin{proposition}\label{PropCodOneSub}
Let $M$ be a simple compact Lie group and let $H$ be a closed connected
subgroup of $M \times M$ which acts hyperpolarly and with cohomogeneity $k \ge
2$ on~$M$. Let $H' \subset H$ be a closed subgroup acting on~$M$ with
cohomogeneity $k + 1$. Then the $H'$-action on $M$ is not polar.
\end{proposition}


\begin{proof}
Assume the $H'$-action on~$M$ is polar. By the results of \cite{hyperpolar}, we
may assume that the $H$-action is a $\sigma$-action or a Hermann action. By
Section~3.2 of \cite{hyperpolar} and Table~5 of~\cite{polar} we know that there
is a point $p \in M$ such that the slice representation of $H_p$ on $\N_p \lp H
\act p \rp$ is irreducible. There are two alternatives for the slice
representation of~$H'_p$ on~$\N_p \lp H' \act p \rp$: Either $\N_p \lp H \act p
\rp = \N_p \lp H' \act p \rp$ and $H'_p$ acts with cohomogeneity $k + 1$ on
$\N_p \lp H \act p \rp$ or $\N_p \lp H \act p \rp$ is a proper submodule of
$\N_p \lp H' \act p \rp$ on which $H'_p$ acts with cohomogeneity~$k$, hence
irreducibly. In the first case a contradiction arises with
Proposition~\ref{PolMinCriteria}~(i). Now consider the second case. Let
$\Sigma'$ be a section of the $H'$-action on $M$ such that $p \in \Sigma'$. It
follows from Corollary~D of~\cite{hl} that $\Sigma'$ is not flat. Since $H'_p$
acts irreducibly on~$\N_p \lp H \act p \rp$, the Weyl group of the slice
representation of~$H'_p$ acts irreducibly on the hyperplane $\T_p \Sigma' \cap
\N_p \lp H \act p \rp$ in $\T_p \Sigma'$. Thus $\Sigma'$ is an irreducible
symmetric space of rank~$k \ge 2$ and dimension~$k + 1$ and we have reached a
contradiction.
\qed\end{proof}

In~\cite{eh2} a list of orbit equivalent subactions of irreducible polar
representations of cohomogeneity~$\ge 2$ is given, we reproduce this list in
Table~\ref{TOrbitEqPolarSubGr} for the convenience of the reader.
%
\begin{table}[h]\rm
  \begin{center}\begin{tabular}{|c|c|c|c|}\hline
    \hl \str $G$ & $K$ & $K'$ & Condition \\\hline\hline
    $\SO(9)$ & $\SO(7) \cdot \SO(2)$ & $\LG_{2} \cdot \SO(2)$ & \\\hline
    $\SO(10)$ & $\SO(8) \cdot \SO(2)$ & $\Spin(7) \cdot \SO(2)$ & \\\hline
    $\SO(11)$ & $\SO(8) \cdot \SO(3)$ & $\Spin(7) \cdot \SO(3)$ & \\\hline
    $\SU(p+q)$ & ${\rm S}(\U(p) \cdot \U(q))$ & $\SU(p) \cdot \SU(q)$ & $p \neq q$ \\\hline
    $\SO(2n)$ & $\U(n)$ & $\SU(n)$ & $n$ odd \\\hline
    $\LE_6$ & $\Spin(10) \cdot \U(1)$ & $\Spin(10)$ & \\\hline
  \end{tabular}

\caption{Orbit equivalent subactions of polar representations.}
\label{TOrbitEqPolarSubGr}
\end{center}
\end{table}
%
Table~\ref{TOrbitEqPolarSubGr} is to be interpreted as follows: If a connected
compact Lie group~$K'$ acts on some finite dimensional Euclidean vector space
by an irreducible polar representation such that the action is non-transitive
on the unit sphere, then either the $K'$-representation is equivalent to an
isotropy representation of a Riemannian symmetric space, or there is a
Riemannian symmetric space $G/K$ and the $K'$-representation is equivalent to
the isotropy representation of~$G/K$ restricted to the subgroup $K' \subset K$
where the triple $(G,K,K')$ is as in Table~\ref{TOrbitEqPolarSubGr}. In the
latter case, the $K$-action and the $K'$-action are orbit equivalent.


\section{Subactions of $\sigma$-actions}
\label{Sigma}


In this section, we do not restrict ourselves to the case of exceptional
groups; we will show for any connected simple compact Lie group of rank greater
than one that for any closed subgroup of~$\Delta^{\sigma} L$ acting polarly
on~$L$ the sections~$\Sigma$ are either flat or $\Sigma = L$.

Let $L$ be a simply connected simple compact Lie group with $\rk L \ge 2$,
equipped with a biinvariant Riemannian metric. Let $\sigma$ be an automorphism
of~$L$ and let $\Delta^{\sigma} L \subset L \times L$ be a subgroup of the form
$$
\Delta^{\sigma} L = \left \{ ( g, \sigma(g) ) \mid g \in L \right \}.
$$
In this section, we will consider the action of $\Delta^{\sigma} L$ and of
closed connected subgroups~$H \subset \Delta^{\sigma} L$ on~$L$. We may assume
that $\sigma$ is either the identity or is induced by a nontrivial automorphism
of the Dynkin diagram of~$L$.

In the first case, the action of $\Delta^{\id_L} L = \Delta L$ is simply the
adjoint action of~$L$. In particular, the identity element of~$L$ is a fixed
point and it follows by Corollary~6.2 of~\cite{hyperpolar}, cf.\ also
\cite{brueck}, Theorem~2.2, that the action of any closed connected subgroup~$H
\subset \Delta L$ on~$L$ is hyperpolar and in fact orbit equivalent to the
$\Delta L$-action since $\rk L \ge 2$. This implies $H = \Delta L$ by
Table~\ref{TOrbitEqPolarSubGr}.

Now let $\sigma$ be an outer automorphism of~$L$ induced by an automorphism of
the Dynkin diagram of~$L$. Then $L$, the order of~$\sigma$, the connected
component of the fixed point group~$L^{\sigma}$ and the cohomogeneity of the
$\Delta^{\sigma} L$-action on~$L$ are as given by Table~\ref{TSigma},
cf.~\cite{hyperpolar}, Theorem~3.4 and \cite{helgason}, Ch.~X, \S5.

\begin{table}[!h]
\begin{center}\begin{tabular}{|c*{5}{|c}|}
\hline
$L$ & $\SU(2n+1)$ & $\SU(2n)$ & $\SO(2n)$ & $\LE_6$ & $\Spin(8)$ \\
\hline
${\rm ord}(\s)$ & $2$ & $2$ & $2$ & $2$ & $3$ \\
\hline
$\lp L^{\sigma} \rp_0$ & $\SO(2n+1)$ & $\mySp(n)$ & $\SO(2n-1)$ & $\LF_4$ & $\LG_2$ \\
\hline
Cohomogeneity & $n$ & $n$ & $n-1$ & $4$ & $2$ \\
\hline
\end{tabular}
\caption{$\sigma$-Actions where $\sigma$ is an outer automorphism.}
\label{TSigma}
\end{center}
\end{table}

As shown in Section~3.2 of~\cite{hyperpolar}, the normal space to the
$\Delta^{\sigma}$-orbit at the identity element of~$L$ is
\begin{equation}\label{SigmaSlice}
\left \{ (X, -X) \mid X \in {\mathfrak l},\; \lp \sigma_* \rp_{\e} (X) = X \right \}
\subset \pp,
\end{equation}
where $\lp \sigma_* \rp_{\e} \colon {\mathfrak l} \to {\mathfrak l}$ denotes the differential
of~$\sigma$ at the identity element of~$L$; furthermore, the slice
representation at $\e \in L$ is equivalent to the adjoint representation of the
fixed point group~$L^{\sigma} = \left \{ g \in L \mid \sigma(g) = g \right \}$.
Assume $H \subset \Delta^{\sigma} L$ is a closed connected subgroup acting
polarly on~$L$.

First assume $\rk(L^{\sigma}) \ge 2$. Consider the isotropy subgroup~$H_{\e}$
of the $H$-action at~$\e$. Its slice representation contains~(\ref{SigmaSlice})
as a submodule. Since $\rk(L^{\sigma}) > 1$, it follows from
Proposition~\ref{PolMinCriteria}~(i) that either the action of $H_{\e}$
on~(\ref{SigmaSlice}) is orbit equivalent to the action of~$\lp \Delta^{\sigma}
L \rp_{\e}$ or the subspace~(\ref{SigmaSlice}) of $\pp = \T_{\e}L$ is tangent
to a section through~$\e$, contradicting Theorem~\ref{ThProductOfSpheres},
since~(\ref{SigmaSlice}) is an irreducible Lie triple system of higher rank.
Thus it follows that the actions of $H_{\e}$ and~$\lp \Delta^{\sigma} L
\rp_{\e}$ on~(\ref{SigmaSlice}) are orbit equivalent. Since the slice
representation of~$\lp \Delta^{\sigma} L \rp_{\e} = \Delta L^{\sigma} \cong
L^{\sigma}$ is equivalent to the adjoint representation of~$L^{\sigma}$, one
can see from Table~\ref{TOrbitEqPolarSubGr} that $\lp H_{\e} \rp_0 = \lp \lp
\Delta^{\sigma} L \rp_{\e}\rp_0$. Since this argument can be applied to all
conjugate subgroups $g H g^{-1} \subset \Delta^{\sigma} L$, $g \in
\Delta^{\sigma} L$, it follows that $H$ contains all conjugates of $\lp \lp
\Delta^{\sigma} L \rp_{\e}\rp_0$. Since $\Delta^{\sigma} L$ is a simple
connected Lie group, this shows that $H = \Delta^{\sigma} L$.

If $\rk(L^{\sigma}) = 1$ then $L \cong \SU(3)$ and $\sigma \colon L \to L$ is
given by complex conjugation. Consider the action of~$\Delta^{\sigma}
\SUxU{2}{1}$ on~$\SU(3)$; it has a slice representation where a one-dimensional
isotropy group acts nontrivially on two two-dimensional submodules; this
representation is non-polar. Any other subgroups of~$\SU(3)$ are of
dimension~$\le 3$; however it is easy to see that $\SU(3)$ does not contain any
totally geodesic subspaces of dimension~$\ge 5$ locally isometric to a product
of spaces of constant curvature. We have shown the following.


\begin{proposition}\label{PropSubSigma}
Let $L$ be simple compact connected Lie group of rank greater than one and let
$\sigma$ be an automorphism of~$L$. Assume the action of a closed connected
non-trivial subgroup~$H$ of $\Delta^{\sigma} L$ on~$L$ is polar. Then $H =
\Delta^{\sigma} L$.
\end{proposition}


Note that the case $\rk L = 1$ is excluded in Proposition~\ref{PropSubSigma}
since the one-dimensional subgroup $\Delta \SUxU{1}{1} \subset \SU(2) \times
\SU(2)$ acts polarly on~$\SU(2) = {\rm S}^3$ with non-flat sections.


\section{Subactions of Hermann actions}
\label{Hermann}


From now on assume that $L$ is an exceptional simple compact Lie group $\LE_6$,
$\LE_7$, $\LE_8$, $\LF_4$, or $\LG_2$. Subgroups in groups of the form
$\Delta^{\sigma}L$ have been treated in Section~\ref{Sigma}.  By
Lemma~\ref{LmMaxSubgr}, we may assume $H'$ is a closed connected subgroup of $H
= H_1 \times H_2 \subset L \times L$, where $H_i \subset L$ are maximal closed
connected subgroups.

We start with subactions of Hermann actions, i.e.\ the special case where both
subgroups $H_1, H_2 \subset L$ are symmetric. All combinations where $H_1$ and
$H_2$ are not conjugate are given by Table~\ref{THermannExc},
cf.~\cite{hyperpolar}. (By $\SO'(2n)$ we denote the image of a half-spin
representation of $\Spin(2n)$.)

\begin{table}[!h]\rm
\begin{center}\begin{tabular}{*{5}{|c}|}
\hline \str
Action&$H_1$&$L$&$H_2$&Coh.\\
\hline\hline
E\,I-II & $\mySp(4)/\{\pm1\}$&$\LE_6$&$\SU(6){\cdot}\mySp(1)$&$4$\\
\hline \hl
E\,I-III& $\mySp(4)/\{\pm1\}$&$\LE_6$&$\Spin(10){\cdot}\U(1)$&$2$\\
\hline \hl
E\,I-IV & $\mySp(4)/\{\pm1\}$&$\LE_6$&$\LF_4$&$2$\\ \hline \hl
E\,II-III & $\SU(6){\cdot}\mySp(1)$&$\LE_6$&$\Spin(10){\cdot}\U(1)$&$2$\\
\hline \hl
E\,II-IV & $\SU(6){\cdot}\mySp(1)$&$\LE_6$&$\LF_4$&$1$\\ \hline \hl
E\,III-IV & $\Spin(10){\cdot}\U(1)$&$\LE_6$&$\LF_4$&$1$\\ \hline \hline \hl
E\,V-VI & $\SU(8)/\{\pm1\}$&$\LE_7$&$\SO'(12){\cdot}\mySp(1)$&$4$\\
\hline \hl
E\,V-VII & $\SU(8)/\{\pm1\}$&$\LE_7$&$\LE_6{\cdot}\U(1)$&$3$\\
\hline \hl
E\,VI-VII & $\SO'(12){\cdot}\mySp(1)$&$\LE_7$&$\LE_6{\cdot}\U(1)$&$2$\\
\hline \hline\hl
E\,VIII-IX & $\SO'(16)$&$\LE_8$&$\LE_7{\cdot}\mySp(1)$&$4$\\\hline \hline \hl
F\,I-II & $\mySp(3){\cdot}\mySp(1)$&$\LF_4$&$\Spin(9)$&$1$\\ \hline
\end{tabular}
\caption{Hermann actions on exceptional groups.} \label{THermannExc}
\end{center}
\end{table}

We will now restrict our attention to subactions of Hermann actions whose
cohomogeneity is greater than one. The actions E\,II-IV, E\,III-IV F\,I-II and
F\,II-II are of cohomogeneity one and their subactions will be treated below.


\begin{theorem}\label{ThmSubHermann}
Let $L$ be an exceptional simple compact Lie group and let $H_1$, $H_2$ be
connected symmetric subgroups of~$L$. Assume that the action of $H$ on $L$ is
of cohomogeneity~$\ge 2$. Let $H'$ be a closed connected nontrivial subgroup of
$H = H_1 \times H_2$. Then the action of~$H'$ on~$L$ is polar if and only if it
is orbit equivalent to the $H$-action. Furthermore, there are orbit equivalent
subactions only for the Hermann actions E\,II-III, E\,III-III and E\,VI-VII.
\end{theorem}


\begin{proof}
First we will consider the special case where $H_1$ and $H_2$ are conjugate; we
may assume $H_1 = H_2$. Then the isotropy subgroup of the $H$-action at the
identity element $\e \in L$ is $H \cap \Delta L = \Delta H_1 = \Delta H_2$ and
the slice representation on $\N_{\e} \lp H \act \e \rp$ is equivalent to the
isotropy representation of the exceptional symmetric space $L / H_1 = L / H_2$.
Assume $H' \subset H$ is a closed connected subgroup acting polarly on~$L$. The
isotropy subgroup of the $H'$-action at the identity element of~$L$ is $H'_{\e}
= H' \cap \Delta L$ and the slice representation of the $H'$-action contains
the normal space $\N_{\e} \lp H \act \e \rp$ to the $H$-orbit through~$\e$ as a
submodule. Since the slice representation of~$H_{\e}$ on  $\N_{\e} \lp H \act
\e \rp$ is irreducible of cohomogeneity~$\ge 2$, it is polarity minimal by
Proposition~\ref{PolMinCriteria}~(i). Thus the action of the group~$H'_{\e}$ on
$\N_{\e} \lp H \act \e \rp$ is either orbit equivalent to the $H_{\e}$-action
or has finite orbits. In the latter case, there arises a contradiction with
Theorem~\ref{ThProductOfSpheres}, since $\N_{\e} \lp H \act \e \rp$ is an
irreducible Lie triple system of higher rank in~$\p$.

Assume the $H'_{\e}$-action on $\N_{\e} \lp H \act \e \rp$ is orbit equivalent
to the $H_{\e}$-action. It then follows from the contents of
Table~\ref{TOrbitEqPolarSubGr} that $H'_{\e}$ contains the connected component
of~$H_{\e} = \Delta H_1$, except possibly in case $L = \LE_6$, $H_1 = H_2 =
\Spin(10) \cdot \U(1)$, where it follows only that $H'_{\e}$ contains a factor
isomorphic to $\Spin(10)$.

Assume for the moment that $H'_{\e}$ contains all of $\Delta H_1$. Let
$(h_1,h_2) \in H_1 \times H_2$ and consider the conjugate subgroup $(h_1,h_2)
\cdot H' \cdot (h_1,h_2)^{-1} \subseteq H$, which also acts polarly on~$L$. It
follows from the above argument that also $(h_1,h_2) \cdot H' \cdot
(h_1,h_2)^{-1}$ contains $\Delta H_1$. This shows that $H'$ contains all
conjugates of~$\Delta H_1$ in~$H$, i.e. $H'$ contains the smallest closed
normal subgroup of~$H$ containing~$\Delta H_1$. In case $H$ is semisimple it
follows that $H' = H$.

Now consider the case where $H$ is not semisimple. There are two cases,
corresponding to the symmetric spaces E\,III and E\,VII.

In case $L = \LE_6$ and $H_1 = \Spin(10) \cdot \U(1)$ it follows from the above
argument that $H'$ contains $\Spin(10) \times \Spin(10)$. Hence $H' \cong \lp
\Spin(10) \times \Spin(10) \rp \cdot Q$, where $Q$ is a closed subgroup
of~$\U(1) \times \U(1)$. Consider the $H'$-orbit $H' \act \e$ through~$\e$; it
is a closed submanifold of codimension~$\le 1$ in the orbit $H \act \e$, which
coincides with the subgroup  $\Spin(10) \cdot \U(1) \subset \LE_6$ and it
follows from Proposition~\ref{PropCodOneSub} that $H' \act \e = H \act \e$. It
follows that either $Q = \U(1) \times \U(1)$ or $Q$ is any one-dimensional
closed subgroup of $\U(1) \times \U(1)$ except the diagonal~$\Delta \U(1)$. In
case $L = \LE_7$ and $H_1 = \LE_6 \cdot \U(1)$, it follows from the above
argument that $H'$ contains $\lp \LE_6 \times \LE_6 \rp \cdot \Delta \U(1)$,
however, this group acts non-polarly on~$\LE_7$ by
Proposition~\ref{PropCodOneSub}. It follows that $H' = H$ in this case.

Now consider the case where $H_1$ and $H_2$ are not conjugate.
%
%
\begin{table}[h]
\begin{center}
\begin{tabular}{|l|c|c|}\hline
\hl \str Action & Slice representation & Kernel \\
\hline\hline \str
E\,I-II    & $\LF_4/\mySp(3)\cdot\mySp(1)$ & $$ \\ \hline \hl
E\,I-III   & $\mySp(4)/\mySp(2)\cdot\mySp(2)$ & $$ \\ \hline \hl
E\,I-IV    & $\SU(6)/\mySp(3)$ & $\mySp(1)$ \\\hline \hl
E\,II-III  &  $\SU(6)/{\rm S}(\U(2)\cdot\U(4))$ & $\mySp(1)$ \\
\hline \hl

E\,II-IV   &  $\mySp(4)/\mySp(3)\cdot\mySp(1)$ & $$ \\ \hline \hl

E\,III-IV  & $\LF_4/\Spin(9)$ & $$ \\ \hline \hline \hl

E\,V-VI    & $\SU(8)/\SUxU{4}{4}$ &  \\ \hline \hl
E\,V-VII   &  $\SU(8)/\mySp(4)$ & \\ \hline \hl
E\,VI-VII  & $\SU(8)/\SUxU{2}{6}$ &  \\ \hline \hline \hl
E\,VIII-IX & $\SO(16)/\U(8)$ & \\ \hline \hline \hl
F\,I-II    & $\mySp(3)/\mySp(2)\cdot\mySp(1)$ & $\mySp(1)$ \\
\hline
\end{tabular}
\caption{Slice representations of Hermann actions on exceptional groups.}
\label{THermannSliceRep}
\end{center}
\end{table}
%
%
We assume that a closed subgroup~$H'$ of $H = H_1 \times H_2$ acts polarly
on~$L$. In Table~\ref{THermannSliceRep}, an irreducible slice representation
for each $H_1 \times H_2$-action is given, cf.~\cite{polar}, Section~3.1.3. The
entries of Table~\ref{THermannSliceRep} are to be interpreted as follows: The
action is given in the first column by the same notation as in
Table~\ref{THermannExc}; in the second column, a symmetric space whose
isotropy representation is (on the Lie algebra level) equivalent to an
effectivized slice representation is given; the third column states the local
isomorphism type of the kernel of the slice representation (if nontrivial).

By an analogous argument as above it follows that the action of~$H'_{\e}$ on
the invariant subspace $\N_{\e} \lp H \act \e \rp$ of the slice representation
is orbit equivalent to the $H_{\e}$-action on~$\N_{\e} \lp H \act \e \rp$,
since we only consider actions of cohomogeneity~$\ge 2$ here.

We start with those Hermann actions where the slice representation restricted
to the connected component of the isotropy group does not have orbit equivalent
proper subgroups. Comparison of Tables~\ref{TOrbitEqPolarSubGr}
and~\ref{THermannSliceRep} shows that this is the case for the actions E\,I-II,
E\,I-III, E\,I-IV, E\,V-VI, E\,V-VII and E\,VIII-IX. We can read off the
isomorphism type of the connected component of the group $H_{\e} = \Delta \lp
H_1 \cap H_2 \rp \subset H$ from Table~\ref{THermannSliceRep}. As above, it
follows that $H'$ contains all conjugates of~$H_{\e}$ in~$H$. In case of the
actions E\,I-II, E\,I-IV, E\,V-VI and E\,VIII-IX this shows that $H = H'$.

In case of the Hermann action~E\,I-III, it follows that $H'$~contains the
subgroup $\lp \SU(4) / \lbc \pm 1 \rbc \rp \times \Spin(10)$ of~$H$. However,
this group acts non-polarly on~$\LE_6$ by Proposition~\ref{PropCodOneSub} and
Theorem~1 of~\cite{polar}, thus $H' = H$.

Similarly, for the action~E\,V-VII, it follows from an analogous argument only
that $H'$ contains the subgroup~$\lp \SU(8) / \lbc \pm 1 \rbc \rp \times \LE_6$
of~$H$ and it follows again from Proposition~\ref{PropCodOneSub} and Theorem~1
of~\cite{polar} that $H' = H$.

For the action~E\,II-III, it follows from the proof of Theorem~7.3
in~\cite{polar} that $H'$ contains the subgroup $\lp \SU(6) \cdot \mySp(1) \rp
\times \Spin(10)$, whose action on~$\LE_6$ is orbit equivalent to the
$H$-action, cf.\ Theorem~2 of~\cite{polar}.

In case of the action~E\,VI-VII we obtain that $H'_{\e}$ contains a
subgroup~$\mySp(1) \cdot \SU(6)$. This shows that $H'$ contains a subgroup $\lp
\SO'(12){\cdot}\mySp(1) \rp \times \LE_6$, whose action on~$\LE_7$ is orbit
equivalent to the $H$-action, cf.\ Theorem~2 of~\cite{polar}.
\qed\end{proof}


\begin{theorem}\label{ThmSubCohOne}
Let $L$ be an exceptional simple compact Lie group and let $H_1$, $H_2$ be
connected subgroups of~$L$ such that $H = H_1 \times H_2$ acts with
cohomogeneity one on~$L$. Let $H'$ be a closed connected nontrivial subgroup of
$H$. Then the action of~$H'$ on~$L$ is either non-polar or orbit equivalent to
the $H$-action.
\end{theorem}


The proof of Theorem~\ref{ThmSubHermann} does not work for actions of
cohomogeneity one, since their slice representations are also of cohomogeneity
one and Proposition~\ref{PolMinCriteria}~(i) cannot be applied.


\begin{proof}
By the results of~\cite{hyperpolar}, there are the following cohomogeneity one
Hermann actions on the simple compact exceptional Lie groups: E\,II-IV,
E\,III-IV, F\,II-II and F\,I-II. We will now treat their subactions case by
case. Here we use the classification of maximal connected subgroups in compact
Lie groups, cf.~\cite{hyperpolar}, Section~2.1.

\paragraph{\bf F\,II-II} Consider subgroups~$H'$ of~$\Spin(9) \times \Spin(9)$
on~$\LF_4$. By Proposition~\ref{PropSubSigma}, we may assume that $H'$ is not
contained in~$\Delta \Spin(9)$. By Lemma~\ref{LmNoLift} and
Lemma~\ref{LmDimSection} we may assume that $H' \subseteq H_1' \times H_2'$,
where $H_1'$ and $H_2'$ are one of the following
\begin{equation}\label{Spin9max1}
\Spin(8),\quad \Spin(7) \cdot \SO(2),\quad \Spin(6) \cdot \Spin(3),\quad
\Spin(5) \cdot \Spin(4).
\end{equation}
First assume $H_1' = H_2'$. Since in this case the slice representation at the
identity element of~$\LF_4$ is equivalent to the isotropy representation of the
homogeneous space $\LF_4 / H_1'$, it follows from Theorem~2 of~\cite{kp} that
the $H'$-action is non-polar.

The remaining actions not excluded by these arguments can be seen to have
non-polar slice representations.

\paragraph{\bf F\,II-I} Assume $H'$ is closed connected subgroup of $H = H_1 \times
H_2 = \Spin(9) \times \lp \mySp(3) \cdot \mySp(1) \rp$ acting polarly on~$\LF_4$.
By Lemma~\ref{LmNoLift} we know that $H'$ does not contain $H_1 \times \{\e\}$
or $\{\e\} \times H_2$. The maximal dimension of a proper closed subgroup in
$H_1$ or $H_2$ is $28$ and $22$, respectively. Hence it follows from
Lemma~\ref{LmDimSection} that $H'$ is contained in a subgroup $H_1' \times
H_2'$ where $H_1'$ is a maximal connected subgroup of $\Spin(9)$ of
dimension~$\ge 18$ and $H_2'$ is a maximal connected subgroup of $\mySp(3) \cdot
\mySp(1)$ of dimension~$\ge 12$. The maximal connected subgroups in~$H_1$ of
dimension~$\ge 18$ are:
\begin{equation}\label{Spin9max2}
\Spin(8),\quad \Spin(7) \cdot \SO(2),\quad \Spin(6) \cdot \Spin(3).
\end{equation}
The maximal connected subgroups of~$H_2$ of dimension~$\ge 12$ are:
\begin{equation}\label{Sp3Sp1max}
\mySp(3) \cdot \U(1), \quad \lp \mySp(2) \cdot \mySp(1) \rp \cdot \mySp(1),\quad \U(3)
\cdot \mySp(1).
\end{equation}
We first consider the action of $\Spin(8) \times \lp \mySp(3) \cdot \mySp(1) \rp$
on~$\LF_4$; it has an isotropy subgroup whose connected component is isomorphic
to $\mySp(2) \cdot \mySp(1)$ and whose slice representation is $\lp {\mathbb{H}}^2
\otimes_{{\mathbb{H}}} {\mathbb{H}}^1 \rp \oplus \R^5$, see \cite{polar}, Section~12, p.~479. This
representation is non-polar~\cite{bergmann} and it can be verified that it is
polarity minimal by looking at the closed subgroups of~$\mySp(2) \cdot \mySp(1)$.
Thus the action on~$\LF_4$ is polarity minimal by
Lemma~\ref{LmPolMinHered}~(ii).

Now consider the action of $\lp \Spin(7) \cdot \SO(2) \rp \times \lp \mySp(3)
\cdot \U(1) \rp$ on~$\LF_4$. An explicit calculation shows that there is an
isotropy subgroup locally isomorphic to $\Spin(5) \cdot \U(1) \cdot \U(1)$
whose non-polar slice representation splits as a direct sum $\R^8 \oplus (\R^5
\otimes \R^2) \oplus \R^2$. This representation is non-polar~\cite{bergmann}.
Using Table~\ref{TOrbitEqPolarSubGr} and Proposition~\ref{PolMinCriteria}~(i)
we see that an $18$-dimensional submodule is polarity minimal.

The action of $\lp \Spin(6) \cdot \Spin(3) \rp \times \lp \mySp(3) \cdot \U(1)
\rp$ has a non-polar slice representation. Its subactions are ruled out by a
dimension count. All other combinations of the groups in (\ref{Spin9max2}) and
(\ref{Sp3Sp1max}) result in actions which are non-polar by
Lemma~\ref{LmDimSection}.

\paragraph{\bf E\,II-IV} Assume $H'$ is a closed subgroup of~$H =
\lp\SU(6)\cdot \mySp(1)\rp \times \LF_4$ acting polarly on~$\LE_6$. By
Lemma~\ref{LmDimSection}, we have $\dim H' \ge 60$. By Lemma~\ref{LmNoLift} it
follows that $H'$ is contained in a subgroup $H_1' \times H_2'$, where $H_1'
\subset \SU(6) \cdot \mySp(1)$ and $H_2' \subset \LF_4$ are maximal connected
subgroups. Since a maximal subgroup of maximal dimension in $\SU(6) \cdot
\mySp(1)$ is $\SU(6) \cdot \U(1)$, we may assume that $H_2' = \Spin(9)$ or
$\mySp(3) \cdot \mySp(1)$, cf.\ Table~\ref{TMaxF4}. The remaining possibilities for
$H_1' \times H_2'$ are
\begin{align*}
&\lp \SUxU{5}{1} \cdot \mySp(1) \rp \times \Spin(9),\quad \lp \mySp(3) \cdot \mySp(1)
\rp \times \Spin(9),\\ &\lp \SU(6) \cdot \U(1) \rp \times \Spin(9),\quad \lp
\SU(6) \cdot \U(1) \rp \times \lp \mySp(3) \cdot \mySp(1) \rp.
\end{align*}
We first consider the cases where $H_2' = \Spin(9)$. From
Table~\ref{THermannSliceRep} we see that there is an isotropy group of the
action E\,II-IV whose connected component is $\mySp(3) \cdot \mySp(1)$. First we
determine the intersection of this group with $\Spin(9)$, from the last entry
of Table~\ref{THermannSliceRep} we see that, possibly after conjugation,
$\mySp(3) \cdot \mySp(1) \cap \Spin(9) = \mySp(2) \cdot \mySp(1) \cdot \mySp(1)$. The
slice representation of this isotropy group contains two submodules equivalent
to ${\mathbb{H}}^2 \otimes_{{\mathbb{H}}} {\mathbb{H}}^1$ and is thus non-polar by \cite{hyperpolar},
Lemma~2.9 and polarity minimal by Proposition~\ref{PolMinCriteria}~(ii). Thus
the $H_1' \times H_2'$-action on~$\LE_6$ is polarity minimal by
Theorem~\ref{ThProductOfSpheres}.

In case of the last group we see from \cite{dynkin1}, Table~25, p.~200, that
there is only one conjugacy class of a subgroup of type~$\LC_3$ in $\LE_6$ and
this has a $3$-dimensional centralizer. It follows that there is an isotropy
subgroup $\mySp(3) \cdot \U(1)$, whose $40$-dimensional slice representation is
the restriction of the isotropy representation of $\LE_6 / (\SU(6) \cdot
\mySp(1))$. This representation is non-polar by
Proposition~\ref{PolMinCriteria}~(i) and Table~\ref{TOrbitEqPolarSubGr}.
Subactions are non-polar by Lemma~\ref{LmDimSection}.

\paragraph{\bf E\,III-IV} By Lemma~\ref{LmNoLift}, we may assume that any compact
subgroup of~$(\Spin(10) \cdot \U(1)) \times \LF_4$ acting polarly with non-flat
sections on~$\LE_6$ is contained in $H_1' \times H_2'$ where $H_1'$ is a
maximal connected subgroup of~$\Spin(10) \cdot \U(1)$ and $H_2'$ is a maximal
connected subgroup of~$\LF_4$. It follows from Lemma~\ref{LmDimSection} that we
may assume $H_1'$ is one of
\begin{align*}
&\Spin(10),\quad \Spin(9)\cdot\U(1),\quad \Spin(8)\cdot\SO(2)\cdot\U(1),\\&\qquad\;\;
\Spin(7)\cdot\Spin(3)\cdot\U(1),\quad
\U(5)\cdot\U(1)
\end{align*}
and that $H_2'$ is one of
\begin{align*}
\Spin(9),\quad \mySp(3) \cdot \mySp(1),\quad \LG_2^1 \cdot \LA_1^8,\quad \SU(3)
\cdot \SU(3),
\end{align*}
cf.\ Table~\ref{TMaxF4}. First assume $H_2' = \Spin(9)$. The group $\Spin(9)$
also occurs as an isotropy group of the $\LF_4$-action on $\LE_6 / \lp
\Spin(10) \cdot \U(1) \rp$, cf.\ Table~\ref{THermannSliceRep}, and we see that
any action of a group $H_1' \times \Spin(9)$ has a slice representation with
two equivalent submodules. Such a representation is non polar by
\cite{hyperpolar}, Lemma~2.9 and polarity minimal by
Lemma~\ref{PolMinCriteria}~(ii). Since the sum of these two submodules is
$32$-dimensional, the $H_1' \times H_2'$-action is non-polar and polarity
minimal, cf.~Lemma~\ref{LmPolMinHered}. This argument shows in particular that
we may assume $\dim H_2' \le 24$. Then the remaining possibilities for $H_1'
\times H_2'$ are
\begin{align*}
& \Spin(10) \times \lp \mySp(3) \cdot \mySp(1) \rp,\quad \lp \Spin(9)\cdot\U(1) \rp
\times \lp \mySp(3) \cdot \mySp(1) \rp, \\
& \Spin(9) \times \lp \mySp(3) \cdot \mySp(1) \rp, \quad \Spin(10) \times \lp \LG_2^1 \cdot
\LA_1^8 \rp
\end{align*}
In case of the first three actions, an isotropy group is $\mySp(2) \cdot \mySp(1)
\cdot \mySp(1)$. We can read off from Table~\ref{THermannSliceRep} that in both
cases the $\mySp(2)$-factor acts nontrivially on at least two factors of the
slice representation, which is hence non-polar~\cite{bergmann}. Subactions can
be excluded by Lemma~\ref{LmDimSection}.

It remains to study the action of $\Spin(10) \times (\LG_2^1 \cdot \LA_1^8)$
on~$\LE_6$; however, by Table~39 of~\cite{dynkin1}, the subgroup $\LG_2^1 \cdot
\LA_1^8$ of $\LF_4$ is contained in the subgroup $\LG_2^1 \cdot \LA_2^{2''}$
of~$\LE_6$. Hence this action is a subaction of {\bf\ref{GrpE61}--\ref{GrpE66}}
which will be shown to be non-polar and polarity minimal in
Section~\ref{PolarE6}.
\qed\end{proof}


\section{Actions on $\LG_2$}
\label{PolarG2}


In this section, we will study those isometric actions on~$\LG_2$ which are
neither subactions of the $\Delta \LG_2$-action nor of the $\SO(4) \times
\SO(4)$-action.

\subsection{Subgroups of $\LG_2$}
\label{SubGrG2}

\begin{proposition}\label{PropSubgroupsOfG2}
All conjugacy classes of closed connected nonabelian subgroups of $\LG_2$ and
their inclusion relations are given by Table~\ref{TSubgroupsOfG2}.
\end{proposition}


\paragraph{\bf Remark.} In Table~\ref{TSubgroupsOfG2}, a tilde is used to distinguish
between nonconjugate isomorphic subgroups; e.g.\ the groups denoted by $\SU(2)$
and $\widetilde{\SU(2)}$ correspond to the subgroups denoted by $\LA_1$ and
$\widetilde{\LA_1}$, respectively, in~\cite{dynkin1}. By the upper indices, the
Dynkin index of subgroups is given. Two subgroups $H_1$, $H_2$ are connected by
a line if and only if there is an element $g \in \LG_2$ such that an inclusion
relation holds between $H_1$ and $g\,H_2\,g^{-1}$.


\begin{table}[!h]
\begin{center}
\begin{picture}(130,100)
\put(0,40){
\begin{tabular}{cccccc}
 & & $\LG_2$ & & & \\
 & & & & & \\
 $\SU(3)$ & & $\SO(4)$ & & & $\LA_1^{28}$ \\
 & & & & & \\
 $\U(2)$ & & $\LA_1^4$ & &  $\widetilde{\U(2)}$ & \\
 & & & & & \\
 $\SU(2)$ & & & &  $\widetilde{\SU(2)}$ & \\
 & & & & & \\
\end{tabular}}
\put(20,16){\line(0,1){14}} \put(20,44){\line(0,1){13}}
\put(70,44){\line(0,1){14}} \put(70,72){\line(0,1){10}}
\put(122,22){\line(0,1){8}} \put(112,49){\line(-4,1){32}}
\put(140,71){\line(-6,1){55}} \put(62,45){\line(-3,1){35}}
\put(24,72){\line(4,1){32}} \put(27,45){\line(3,1){35}}
\end{picture}
\caption{Conjugacy classes of nonabelian connected subgroups in~$\LG_2$.}
\label{TSubgroupsOfG2}
\end{center}
\end{table}


\begin{proof}
It is straightforward to prove the proposition using the results
of~\cite{dynkin1}. The conjugacy classes of three dimensional (hence simple)
connected subgroups of~$\LG_2$ are given in Table~16 of~\cite{dynkin1}; there
are four classes, distinguished by their Dynkin indices, which are $1$, $3$,
$4$ and $28$. These subgroups are denoted in~\cite{dynkin1} by $\LA_1$, $\tilde
\LA_1$, $\LA_1^4$ and $\LA_1^{28}$, respectively. There are two conjugacy
classes of maximal regular connected subgroups, $\SU(3) = \LA_2^1$ and $\SO(4)
= \LA_1 \cdot \tilde \LA_1$; moreover, there is only one conjugacy class of
connected subgroups not contained in a proper regular subgroup and this
is~$\LA_1^{28}$. The maximal connected subgroups of $\SO(4) = \LA_1 \cdot
\tilde \LA_1$ are $\SO(3)$ and the two subgroups which we denote by $\U(2)$ and
$\widetilde{\U(2)}$, containing the two simple factors $\LA_1$ and $\tilde
\LA_1$ of~$\SO(4)$, which are not conjugate in $\LG_2$, since they have
different Dynkin indices. It follows from Table~16 of~\cite{dynkin1} that
$\SO(3)$ corresponds to the group $\LA_1^4$. The maximal connected subgroups
of~$\SU(3)$ are $\SUxU{2}{1}$ and $\SO(3)$. Since the first group has Dynkin
index one as a subgroup of~$\LG_2$, it follows that it corresponds to the
subgroup denoted by $\U(2)$ in Table~\ref{TSubgroupsOfG2}, the second group
obviously corresponds to~$\LA_1^4$.
\qed\end{proof}


It follows from Table~\ref{TSubgroupsOfG2} that all connected proper subgroups
of~$\LG_2$ except $\SU(3)$ and $\LA_1^{28}$ are contained in $\SO(4)$ after
conjugation with a suitable element from~$\LG_2$. Since the $\SO(4) \times
\SO(4)$-action and the $\Delta \LG_2$-action were already shown to be polarity
minimal, it suffices to consider the actions of subgroups~$H \subseteq H_1
\times H_2$ where at least one of the closed connected subgroups $H_1, H_2
\subsetneq \LG_2$ is conjugate to either $\SU(3)$ or $\LA_1^{28}$. Let $\pi_i
\colon (g_1, g_2) \mapsto g_i$ for $i = 1,2$ be the canonical projections
$\LG_2 \times \LG_2 \to \LG_2$.

Let us first consider the case where at least one of the factors $\pi_1(H)$ or
$\pi_2(H)$ is conjugate to~$\SU(3)$. We may assume w.l.o.g.\ $\pi_2(H) =
\SU(3)$. Now if $\pi_1(H)$ is conjugate to $\SU(3)$, too, then $H$ is conjugate
to either $\SU(3) \times \SU(3)$ or $\Delta \SU(3)$; in the first case, the
action is a well known cohomogeneity one action, in the latter case, the
$H$-action on~$\LG_2$ is non-polar by Proposition~\ref{PropSubSigma}, since it
is not orbit equivalent to the $\Delta \LG_2$-action. If $\pi_1(H)$ is not
conjugate to $\SU(3)$, then it follows that $H$ is of the form $H_1 \times
H_2$, where $H_1 = \ker \pi_2|_H$ and $H_2 = \ker \pi_1|_H \cong \SU(3)$. We
will consider this case in Subsection~\ref{SU3Act}. It remains the case where
one of the factors $\pi_1(H)$ or $\pi_2(H)$ is conjugate to~$\LA_1^{28}$, which
we will treat in Subsection~\ref{SubA128}.

\subsection{Actions of $H_1 \times \SU(3)$ on $\LG_2$}
\label{SU3Act}

\begin{lemma}\label{LmG2SU3}
Let $H_1 \subset \LG_2$ be a closed subgroup. If the cohomogeneity of the
action of $H_1 \times \SU(3)$ on~$\LG_2$ is greater than one then the action is
not polar.
\end{lemma}


\begin{proof}
Assume the action of $H_1 \times H_2$ on~$\LG_2$ is polar, where $H_2 =
\SU(3)$. Let $\h_i$ be the Lie algebra of~$H_i$ and let $\m_i$ be the
orthogonal complement of~$\h_i$ in~$\Lg_2$ for $i = 1,2$. We may assume that
the identity element $\e \in \LG_2$ lies in a principal orbit. Using
Proposition~\ref{PropPolCrit}~(ii), it follows in particular that $p\lp[\nu,
\nu]\rp = 0$, where $\nu = \m_1 \cap \m_2 = \N_{\e} \lp H_1 \times H_2 \rp \act
\e.$ and where $p$ denotes the orthogonal projection $\Lg_2 \to \h_2$. We want
to study the $\R$-bilinear map $\beta \colon \m_2 \times \m_2 \to \h_2$ given
by $\lp x, y \rp \mapsto p([x,y])$.

To describe this map explicitly, consider the antisymmetric bilinear map $F
\colon \C^3 \times \C^3 \to {\mathfrak u}(3)$, $(x,y) \mapsto x\bar{y}^t-y\bar{x}^t$,
which maps a pair of vectors in~$\C^3$ to a skew-hermitian matrix. Define $F_0
\colon \C^3 \times \C^3 \to \su(3)$ by $F_0(x,y) := F(x,y) - \frac 1 3 \tr \lp
F(x,y) \rp I$, where $I$ denotes the $3 \times 3$ identity matrix. Let $A \in
\SU(3)$; then $F_0( Ax, Ay ) = A \cdot F(x,y) \cdot A^{-1}$, i.e.\ $F_0$
defines a non-zero, hence surjective, $\SU(3)$-equivariant map $\Lambda^2 \R^6
\to \su(3)$.

The adjoint representation of~$\LG_2$ restricted to~$H_2$ leaves $\m_2$
invariant and the action of~$H_2$ on~$\m_2$ is equivalent to the natural
$\SU(3)$-representation on~$\C^3 = \R^6$. We may thus identify $\m_2$ with
$\C^3$ by an $\SU(3)$-equivariant $\R$-linear isomorphism.  Since $\Lambda^2
\R^6$, considered as an $\SU(3)$-module, contains only one irreducible summand
equivalent to the adjoint representation of~$\SU(3)$, it follows from Schur's
Lemma that $\beta$ and $F_0$ agree up to an equivariant isomorphism ($\beta$ is
non-zero since otherwise $\m_2$ would be an ideal of~$\Lg_2$).

Now assume $x,y$ are two non-zero elements in $\nu \subseteq \m_2 = \C^3$ such
that $F_0(x,y)=0$. We will show that $x$ and $y$ are linearly dependent
over~$\R$. Since $F_0$ is $\R$-bilinear and $\SU(3)$ acts transitively on the
unit sphere in~$\R^6$, we may assume that $x = \lp 1,0,0 \rp^t \in \C^3$. Then
we have
$$
F_0 (x,y) = \begin{pmatrix}
              \frac23 \lp \bar y_1 - y_1 \rp & \bar y_2 & \bar y_3 \\
              -y_2 & -\frac13 \lp \bar y_1 - y_1 \rp  & 0 \\
              -y_3 & 0 & -\frac13 \lp \bar y_1 - y_1 \rp \\
            \end{pmatrix},
$$
where $y = \lp y_1, y_2, y_3 \rp^t$; this shows that $F_0(x,y) = 0$ if and only
if $x = \alpha\, y$ for some $\alpha \in \R$. Hence $\nu$ is at most
one-dimensional and the cohomogeneity of the action is at most one.
\qed\end{proof}


Assume now $H = H_1 \times \SU(3)$ acts polarly, hence with cohomogeneity one,
on~$\LG_2$. It follows that $\dim( H_1 ) \ge 5$. A glance at
Table~\ref{TSubgroupsOfG2} now shows that $H_1 = \SU(3)$ or $\SO(4)$. The
actions of these groups are of cohomogeneity one~\cite{hyperpolar}.

\subsection{Subactions of $H_1 \times \LA_1^{28}$ on~$\LG_2$}
\label{SubA128}

Now consider the case where one of the factors $\pi_1(H)$ or $\pi_2(H)$ is
conjugate to~$\LA_1^{28}$. Assume $H$ acts polarly on~$\LG_2$. By
Lemma~\ref{LmDimSection}, we have $\dim(H) \ge 8$ and the only possibilities
for $H$ are $\LA_1^{28} \times \SU(3)$ and $\LA_1^{28} \times \SO(4)$. The
first action is non-polar by Lemma~\ref{LmG2SU3}. It remains to study the
action of $H = \SO(4) \times \LA_1^{28}$ on~$\LG_2$. This action is non-polar
by Lemma~\ref{LmNoLift} and the results of~\cite{hyperpolar}. Since any proper
closed subgroup of~$H$ is of dimension~$\le 7$, the action is polarity minimal
by Lemma~\ref{LmDimSection}.


\section{Actions on~$\LF_4$}
\label{PolarF4}


We will now study actions on~$\LF_4$ which are neither subactions of
$\sigma$-actions nor of Hermann actions. By Lemma~\ref{LmMaxSubgr} we may
assume that the group acting polarly is contained in $H_1 \times H_2$, where
$H_i$ are maximal connected subgroups. By Lemma~\ref{LmDimBd} it follows that
$\dim H_i \ge 12$. The conjugacy classes of all such subgroups of~$\LF_4$ are
given by Table~\ref{TMaxF4}, cf.~\cite{dynkin1}. In the case of a symmetric
subgroup the type of the symmetric space is given in the last column. Of course
we do not need to consider groups $H_1 \times H_2$ where $H_1$ and $H_2$ are
both symmetric subgroups, since they have already been considered in
Section~\ref{Hermann}. We will follow the same procedure in
Sections~\ref{PolarE6} and~\ref{PolarE7} for actions on the groups $\LE_6$ and
$\LE_7$.

\begin{table}[h]\rm
\begin{center}\begin{tabular}{|c|c|c|l|}
\hline \str No. & Subgroup & Dimension & Type \\ \hline \hline

 \mylabel{GrpF41} & $\Spin(9)$ & $36$ & F\,II \\ \hline

 \mylabel{GrpF42} & $\mySp(3) \cdot \mySp(1)$ & $24$ & F\,I \\ \hline

 \mylabel{GrpF43} & $\LG_2^1 \cdot \LA_1^8$ & $17$ & \\ \hline

 \mylabel{GrpF44} & $\SU(3) \cdot \SU(3)$ & $16$ & \\ \hline

\end{tabular}
\caption{Maximal connected subgroups of~$\LF_4$ of dimension~$\ge12$.}
\label{TMaxF4}
\end{center}
\end{table}

\paragraph{\bf\ref{GrpF41}--\ref{GrpF43}.}
We determine a slice representation of the action of $H = H_1 \times H_2$,
where $H_1 = \LG_2^1 \cdot \LA_1^8$ and $H_2 = \Spin(9)$, on~$\LF_4$. Consider
the subgroup $\LG_2 \subset \Spin(7) \subset \Spin(9)$. By~\cite{dynkin1},
Table~25, p.~199, there is only one conjugacy class of subalgebras isomorphic
to~$\LG_2$ and it follows that there is an isotropy subgroup containing the
$\LG_2$-factor of the group \ref{GrpF43}. By Table~25 in~\cite{dynkin1}, the
dimension of the normal space is a multiple of~$7$. A dimension count shows
that such an isotropy group is of dimension~$15$ and thus its Lie algebra is
isomorphic to $\Lg_2 \oplus \R$. From the fact that $\LF_4 / \LG_2^1 \cdot
\LA_1^8$ is a strongly isotropy irreducible homogeneous space \cite{wolfIrr},
one can deduce that the $14$-dimensional slice representation is orbit
equivalent to the action of $\LG_2 \times \SO(2)$ on $\R^7 \otimes \R^2$ given
by the tensor product of the $7$-dimensional irreducible $\LG_2$-representation
and a non-trivial $2$-dimensional real representation of $\SO(2)$. Thus the
action is of cohomogeneity two (see Table~\ref{TOrbitEqPolarSubGr}) and
non-polar by Lemma~\ref{LmNoLift}. By Table~\ref{TOrbitEqPolarSubGr} and since
the $14$-dimensional slice representation is polarity minimal, it follows that
any closed subgroup~$H'$ of $H$ acting polarly on~$\LF_4$ must contain $\Delta
\LG_2^1$. An argument similar as in the proof of Theorem~\ref{ThmSubHermann}
now shows that $H'$ contains $\LG_2^1 \times \Spin(9)$. Hence the $H'$-action
is non-polar by Lemma~\ref{LmNoLift}.

\paragraph{\bf\ref{GrpF41}--\ref{GrpF44}.}
The action of $H_2 = \SU(3) \cdot \SU(3)$ on the Cayley plane $\LF_4 /
\Spin(9)$ is polar of cohomogeneity two, see \cite{pth1}, hence with non-flat
sections. It follows from Lemma~\ref{LmNoLift} that the action of $\Spin(9)
\times H_2'$ on~$\LF_4$ is non-polar for all closed subgroups $H_2' \subseteq
\SU(3) \cdot \SU(3)$. Assume a group $H'$ acting polarly on~$\LF_4$ is
contained in $H_1' \times \lp \SU(3) \cdot \SU(3) \rp$ where $H_1'$ is a
maximal connected subgroup of~$\Spin(9)$. By Lemma~\ref{LmDimSection} it
follows that $H_1' = \Spin(8)$. Consider the action of $\Spin(8) \times \lp
\SU(3) \cdot \SU(3) \rp$ on~$\LF_4$. By a calculation as in~\cite{polar},
Remark~10.1, one finds an isotropy subgroup $\SU(3) \cdot T^2 = \U(3) \cdot
\SO(2)$. The $18$-dimensional slice representation, when restricted
to~$\SU(3)$, splits into three times the standard representation on~$\R^6$.
This representation is non-polar and polarity minimal~\cite{bergmann},
\cite{dadok}.

\paragraph{\bf\ref{GrpF42}--\ref{GrpF43} and \ref{GrpF42}--\ref{GrpF44}.}
Let $H_1 = \mySp(3) \cdot \mySp(1)$ and let $H_2 = \LG_2^1 \cdot \LA_1^8$ or
$\SU(3) \cdot \SU(3)$. By the results of~\cite{polar} and Lemma~\ref{LmNoLift},
the action of $H = H_1 \times H_2$ on~$\LF_4$ is non-polar. Since any proper
closed subgroup of~$H$ is of dimension~$\le 39$, it follows from
Lemma~\ref{LmDimSection} that no closed connected nontrivial subgroup of~$H$
acts polarly on~$\LF_4$.


\section{Actions on~$\LE_6$}
\label{PolarE6}


We follow the same procedure as in Section~\ref{PolarF4}. The maximal connected
subgroups of dimension greater or equal to~$15$ are given in
Table~\ref{TMaxE6}.

\begin{table}[h]\rm
\begin{center}\begin{tabular}{|c|c|c|l|}
\hline \str No. & Subgroup & Dimension & Type \\ \hline \hline

 \mylabel{GrpE64} & $\LF_4$ & $52$ & E\,IV \\ \hline

 \mylabel{GrpE61} & $\Spin(10) \cdot \U(1)$ & $46$ & E\,III \\ \hline

 \mylabel{GrpE62} & $\SU(6) \cdot \mySp(1)$ & $38$ & E\,II \\ \hline

 \mylabel{GrpE65} & $\mySp(4) / \{\pm1\}$ & $36$ & E\,I \\ \hline

 \mylabel{GrpE63} & $\SU(3) \cdot \SU(3) \cdot \SU(3)$ & $24$ & \\ \hline

 \mylabel{GrpE66} & $\LG_2^1 \cdot \LA_2^{2''}$ & $22$ & \\ \hline

\end{tabular}
\caption{Maximal connected subgroups of~$\LE_6$ of dimension~$\ge 15$.}
\label{TMaxE6}
\end{center}
\end{table}

\paragraph{\bf\ref{GrpE64}--\ref{GrpE63} and \bf\ref{GrpE64}--\ref{GrpE66}.}
For both actions, a slice representation is computed in~\cite{polar},
Subsection~10.1.  It is shown there that these representations are non-polar
and polarity minimal. Since these slice representations are of dimension $36$
and $21$, respectively, it follows from Lemma~\ref{LmPolMinHered}~(ii) that
both actions on $\LE_6$ are non-polar and polarity minimal.

\paragraph{\bf\ref{GrpE61}--\ref{GrpE63}.}
See Section~\ref{Regular}.

\paragraph{\bf\ref{GrpE61}--\ref{GrpE66}.}
Let $H_1 = \Spin(10) \cdot \U(1)$, $H_2 = \LG_2^1 \cdot \LA_2^{2''}$. Using
Table~25 of~\cite{dynkin1}, pp.~200 and 203, it follows that $\LG_2^1$ is
contained (after conjugation) in~$H_1$. Since the isotropy representation of
the strongly isotropy irreducible space $\LE_6 / (\LG_2^1 \cdot \LA_2^{2''})$
is equivalent to the tensor product of the adjoint representation of~$\SU(3)$
and the $7$-dimensional irreducible representation of $\LG_2$, it follows that
the dimension of the slice representation of $H_1 \cap H_2$ is a multiple
of~$7$. A dimension count now shows that the isotropy group $H_1 \cap H_2$ is
locally isomorphic to $\LG_2 \cdot \SUxU{2}{1}$ and the $28$-dimensional slice
representation is non-polar~\cite{dadok}, \cite{bergmann} and polarity minimal.
Thus the $H_1 \times H_2$-action is non-polar and polarity minimal by
Proposition~\ref{LmPolMinHered}~(ii).

\paragraph{\bf\ref{GrpE62}--\ref{GrpE66} and \ref{GrpE65}--\ref{GrpE63}.}
The actions of the groups $H = \lp \SU(6) \cdot
\mySp(1) \rp \times ( \LG_2^1 \cdot \LA_2^{2''} )$ and $\lp \mySp(4) / \{\pm1\} \rp
\times \lp \SU(3) \cdot \SU(3) \cdot \SU(3) \rp$ are non-polar by
Lemma~\ref{LmNoLift} and the results of~\cite{hyperpolar}. Since both groups
are $60$-dimensional, it follows from Lemma~\ref{LmDimSection} that no closed
connected nontrivial subgroup of these groups acts polarly on~$\LE_6$.

\paragraph{\bf\ref{GrpE62}--\ref{GrpE63}.}
See Section~\ref{Regular}.


\section{Actions on~$\LE_7$} \label{PolarE7}


The maximal connected subgroups of~$\LE_7$ of dimension~$\ge 34$ are given in
Table~\ref{TMaxE7}, cf.~\cite{dynkin1}. The groups \ref{GrpE78}, \ref{GrpE71}
and \ref{GrpE74} are the symmetric subgroups of~$\LE_7$.

\begin{table}[h]\rm
\begin{center}\begin{tabular}{|c|c|c|l|}
\hline \str No. & Subgroup & Dimension & Type \\ \hline \hline

 \mylabel{GrpE78} & $\LE_6 \cdot \U(1)$ & $79$ & E\,VII \\ \hline

 \mylabel{GrpE71} & $\SO'(12) \cdot \mySp(1)$ & $69$ & E\,VI \\ \hline

 \mylabel{GrpE74} & $\SU(8)/\{\pm1\}$ & $63$ & E\,V \\ \hline

 \mylabel{GrpE75} & $\LF_4^1 \cdot \LA_1^{3''}$ & $55$ & \\ \hline

 \mylabel{GrpE72} & $\SU(6) \cdot \SU(3)$ & $43$ & \\ \hline

 \mylabel{GrpE76} & $\LG_2^1 \cdot \LC_3^{1''}$ & $35$ & \\ \hline

\end{tabular}
\caption{Maximal connected subgroups of~$\LE_7$ of dimension~$\ge 34$.}
\label{TMaxE7}
\end{center}
\end{table}

\paragraph{\bf\ref{GrpE78}--\ref{GrpE75} and \bf\ref{GrpE71}--\ref{GrpE75}.}
Let $H_1 = \LF_4^1 \cdot \LA_1^{3''}$. For the cases $H_2 = \LE_6 \cdot \U(1)$
or $H_2 = \SO'(12) \cdot \mySp(1)$, the argument given in~\cite{polar},
Section~10.2 shows that the $H$-action on~$L$ is non-polar and polarity minimal
by Lemma~\ref{LmPolMinHered}~(ii).

\paragraph{\bf\ref{GrpE78}--\ref{GrpE72} and \ref{GrpE71}--\ref{GrpE72}.}
See Section~\ref{Regular}.

\paragraph{\bf\ref{GrpE78}--\ref{GrpE76}.}
It follows from \cite{dynkin1}, Table~25, that the group $\LG_2^1$ is contained
(after conjugation) in $\LE_6 \subset \LE_7$ and that the dimension of a slice
representation of an isotropy subgroup containing $\LG_2^1$ is a multiple
of~$7$. Now consider the subgroup $\LG_2^1 \cdot \LA_2^{2''}$ of $\LE_6$. It
can be read off from Table~25 in~\cite{dynkin1} that this group is contained in
the subgroup $\LG_2^1 \cdot \LC_3^{1''}$ of $\LE_7$. A dimension count shows
that an isotropy group containing $\LG_2^1 \cdot \LC_3^{1''}$ must be locally
isomorphic to $\LG_2 \cdot \U(3)$. The corresponding slice representation is
equivalent to the $42$-dimensional real tensor product of the $7$-dimensional
$\LG_2$-representation and the real $6$-dimensional standard
$\SU(3)$-representation. This slice representation is irreducible and
non-polar~\cite{dadok}, hence polarity minimal by
Proposition~\ref{PolMinCriteria}. We conclude that the action
{\bf\ref{GrpE78}--\ref{GrpE76}} is non-polar and polarity minimal by
Lemma~\ref{LmPolMinHered}~(ii).

\paragraph{\bf\ref{GrpE74}--\ref{GrpE75}.}
The action of $\lp \SU(8)/\{\pm1\} \rp \times H_2'$ on~$\LE_7$, where $H_2'$ is
a closed subgroup of~$\LF_4^1 \cdot \LA_1^{3''}$, is non-polar by
Lemma~\ref{LmNoLift}. The dimension of any other proper closed subgroup in $\lp
\SU(8)/\{\pm1\} \rp \times \LF_4^1 \cdot \LA_1^{3''}$ is less than~$105$.


\section{Actions on~$\LE_8$}
\label{PolarE8}


Since all closed connected subgroups of~$\LE_8$ of dimension~$\ge 90$ are
symmetric~\cite{dynkin1}, it follows from Lemma~\ref{LmDimBd} and
Theorem~\ref{ThmSubHermann} that any polar action on $\LE_8$ is hyperpolar or
has finite orbits.


\section{Regular subgroups of the isometry group}
\label{Regular}


\begin{table}[h]\tiny
\begin{center}\begin{tabular}{*{3}{|c}|}
\hline
\str Action & Isotropy subgroup & Slice representation \\
\hline
\hline
{\mbox{\bf\ref{GrpE61}--\ref{GrpE63}}} &
%
\begin{minipage}{75pt}
\begin{picture}(100,60)
\put(38,17){\circle{3}} \multiput(6,17)(16,0){4}{\circle*{3}}
\multiput(6,17)(16,0){2}{\circle{5}} \multiput(38,33)(0,16){2}{\circle{5}}
\multiput(54,17)(16,0){2}{\circle{5}}
\multiput(8.5,17)(16,0){4}{\line(8,0){11}}
\multiput(38,19.5)(0,16){2}{\line(0,8){11}}
 \put(38,33){\circle*{3}}
\put(4,5){$$} \put(20,5){$$} \put(36,5){$$} \put(52,5){$$} \put(68,5){$$}
\put(44,30){$$}
\end{picture}
\end{minipage} &
$ \begin{picture}(60,20) \multiput(6,2)(16,0){4}{\circle*{3}}
\multiput(6,2)(16,0){4}{\circle{5}} \multiput(8.5,2)(16,0){1}{\line(8,0){11}}
\put(3.4,7){$1$} \put(35.4,7){$1$}
\end{picture}
\oplus  \begin{picture}(60,20) \multiput(6,2)(16,0){4}{\circle*{3}}
\multiput(6,2)(16,0){4}{\circle{5}} \multiput(8.5,2)(16,0){1}{\line(8,0){11}}
\put(19.4,7){$1$} \put(51.4,7){$1$}
\end{picture} $
\\
\hline
{\mbox{\bf\ref{GrpE62}--\ref{GrpE63}}} &
%
\begin{minipage}{78pt}
\begin{picture}(100,60)
\multiput(6,17)(16,0){5}{\circle*{3}} \multiput(6,17)(16,0){2}{\circle{5}}
\multiput(38,33)(0,16){2}{\circle{5}} \multiput(54,17)(16,0){2}{\circle{5}}
\multiput(8.5,17)(16,0){4}{\line(8,0){11}}
\multiput(38,19.5)(0,16){2}{\line(0,8){11}} \put(38,49){\circle*{3}}
\put(4,5){$$} \put(20,5){$$} \put(36,5){$$} \put(52,5){$$} \put(68,5){$$}
\put(44,30){$$}
\end{picture}
\end{minipage} &
$ \begin{picture}(78,20) \multiput(6,2)(16,0){5}{\circle*{3}}
\multiput(6,2)(16,0){5}{\circle{5}} \multiput(8.5,2)(32,0){2}{\line(8,0){11}}
\put(3.4,7){$1$} \put(51.4,7){$1$} \put(67.4,7){$1$}
\end{picture}
$
\\
\hline
{\mbox{\bf\ref{GrpE78}--\ref{GrpE72}}} &
%
\begin{minipage}{108pt}
\begin{picture}(100,45)
\multiput(22,17)(16,0){5}{\circle*{3}} \multiput(6,17)(16,0){4}{\circle{5}}
\multiput(8.5,17)(16,0){6}{\line(8,0){11}}
\multiput(54,19.5)(0,16){1}{\line(0,8){11}}
\multiput(86,17)(16,0){2}{\circle{5}} \put(54,33){\circle{5}}
\put(54,33){\circle*{3}} \put(4,5){$$} \put(20,5){$$} \put(36,5){$$}
\put(52,5){$$} \put(68,5){$$} \put(44,30){$$}
\end{picture}
\end{minipage} &
$ \begin{picture}(78,20) \multiput(6,2)(16,0){5}{\circle*{3}}
\multiput(6,2)(16,0){5}{\circle{5}} \multiput(8.5,2)(16,0){3}{\line(8,0){11}}
\put(3.4,7){$1$} \put(67.4,7){$1$}
\end{picture}
\oplus  \begin{picture}(78,20) \multiput(6,2)(16,0){5}{\circle*{3}}
\multiput(6,2)(16,0){5}{\circle{5}} \multiput(8.5,2)(16,0){3}{\line(8,0){11}}
\put(35.4,7){$1$}
\end{picture} $
\\

\hline
{\mbox{\bf\ref{GrpE71}--\ref{GrpE72}}} &
%
\begin{minipage}{108pt}
\begin{picture}(100,45)
\multiput(6,17)(16,0){5}{\circle*{3}} \multiput(6,17)(16,0){4}{\circle{5}}
\multiput(8.5,17)(16,0){6}{\line(8,0){11}}
\multiput(54,19.5)(0,16){1}{\line(0,8){11}}
\multiput(86,17)(16,0){2}{\circle{5}} \put(54,33){\circle{5}}
\put(54,33){\circle*{3}} \put(102,17){\circle*{3}} \put(4,5){$$} \put(20,5){$$}
\put(36,5){$$} \put(52,5){$$} \put(68,5){$$} \put(44,30){$$}
\end{picture}
\end{minipage}&
$ \begin{picture}(94,20) \multiput(6,2)(16,0){6}{\circle*{3}}
\multiput(6,2)(16,0){6}{\circle{5}} \multiput(8.5,2)(16,0){4}{\line(8,0){11}}
\put(19.4,7){$1$} \put(83.4,7){$1$}
\end{picture}
$
\\
\hline
\end{tabular}
\caption{Isotropy groups and slice representations of certain regular
subgroups.}\label{TRegAct}
\end{center}
\end{table}

For the special case where $H_1,H_2$ are maximal regular subgroups of the
simple compact Lie group~$G$, one can determine a slice representation of the
$H_1 \times H_2$-action on~$G$ by the method described in Remark~10.1
of~\cite{polar}, see also \S3 of~\cite{oniscikBook}, Theorem~16. We will apply
this method now to certain actions on~$\LE_6$ and $\LE_7$.

In the middle column of Table~\ref{TRegAct}, the extended Dynkin diagram of~$G$
is given for each action under consideration. We assume that the intersection
of $H_1$ and $H_2$ contains a fixed maximal torus~$T$ of~$G$. Then the root
systems of $H_1$ and $H_2$ with respect to~$T$ are subsets~$S_1$, $S_2$ of the
root system~$R$ of~$G$. Simple roots of the root systems of $H_1$ or $H_2$ are
shown in Table~\ref{TRegAct} by black nodes~\hbox{{\put(5,3){\circle*{3}}}{\ \
\,}} or circles~\hbox{{\put(5,3){\circle{5}}}{\ \ \;}}, respectively. The
intersection~$S_1 \cap S_2$ of both root system is the root system of the
isotropy group~$H_1 \cap H_2$, its simple roots are shown
as~\hbox{{\put(5,3){\circle*{3}}\put(5,3){\circle{5}}}{\ \ \;}}. The roots in
$R \setminus \lp S_1 \cup S_2 \rp$ are exactly the weights of the slice
representation of~$H_1 \cap H_2$. In the third column of Table~\ref{TRegAct},
highest weights of the irreducible submodules of the slice representation,
restricted to the semisimple part of $H_1 \cap H_2$, and viewed as complex
representations, are given. The slice representations are of dimension $24$,
$36$, $40$ and $60$, respectively; they are non-polar~\cite{dadok},
\cite{bergmann} and polarity minimal by Proposition~\ref{PolMinCriteria}. Thus
all four actions are non-polar and polarity minimal by
Lemma~\ref{LmPolMinHered}~(ii).


\section{Principal isotropy algebras of actions on~$\LG_2$}
\label{G2Isotropy}


In this section, we determine all isometric actions on the compact exceptional
Lie group $\LG_2$ which have principal isotropy subgroups of positive
dimension.


\begin{theorem}\label{ThmPrOrbG2}
Let $H \subseteq \LG_2 \times \LG_2$ be a closed connected subgroup acting
nontransitively on~$\LG_2$. Then the principal isotropy groups of the
$H$-action on~$\LG_2$ are finite except if $H$ is conjugate to $\Delta \LG_2$,
$\SU(3) \times \SU(3)$, $\SU(3) \times \SO(4)$, or $\SO(4) \times \SU(3)$.
\end{theorem}


In particular, it follows from Theorem~\ref{ThmPrOrbG2} and
Table~\ref{TSubgroupsOfG2} that for all integers $0 \le d \le 14$ there is a
closed subgroup $H \subset \LG_2 \times \LG_2$ such that the $H$-action on
$\LG_2$ has principal orbits of dimension~$d$.


\begin{proof}
We will use the same kind of recursion procedure as in the proof of
Theorem~\ref{ThmCohOne} to classify all closed subgroups $H \subset \LG_2
\times \LG_2$ such that $H$ acts nontransitively on $\LG_2$ with principal
isotropy groups of positive dimension. By Lemma~\ref{LmMaxSubgr}, we may assume
that either $H$ is contained in $\Delta \LG_2$ or that $H$ is contained in a
group of the form $H_1 \times H_2$ where $H_i \subset \LG_2$ are maximal
connected subgroups, cf.\ Table~\ref{TSubgroupsOfG2}.

Assume first $H \subseteq \Delta \LG_2$. For the actions of such groups the
identity element $\e \in \LG_2$ is a fixed point and we may consider the slice
representations of these actions at~$\e$, which are equivalent to the adjoint
representation of~$\LG_2$ restricted to~$H$. The action of $\Delta \LG_2$ on
$\LG_2$ is the adjoint action whose principal isotropy groups are the maximal
tori of~$\LG_2$. Now consider the maximal connected subgroups $H' = \Delta
\SU(3)$, $\Delta \SO(4)$ and $\Delta \LA_1^{28}$ of $\Delta \LG_2 \cong \LG_2$,
cf.\ Table~\ref{TSubgroupsOfG2}. In all three cases the adjoint representation
of~$\LG_2$ restricted to~$H'$ is equivalent to the adjoint representation
of~$H'$ plus the isotropy representation of the strongly isotropy irreducible
homogeneous space $\LG_2 / H'$~\cite{wolfIrr}. We see from~\cite{hh} that these
representations have finite principal isotropy groups. Hence all subactions of
these action also have finite principal isotropy groups.

Now consider $H$ contained in groups of the form $H' = H_1 \times H_2$, where
$H_i \in \lbc \SU(3), \SO(4), \LA_1^{28} \rbc$. The action of $\SO(4) \times
\SO(4)$ on~$\LG_2$ has finite principal isotropy groups; therefore, in order to
find all closed subgroups of $H'= H_1 \times H_2$ acting with principal
isotropy groups of positive dimension, we may assume that $H_1 = \SU(3)$ or
$H_1 = \LA_1^{28}$ and that $H$ contains the $H_1$-factor of $H'$. We will
treat the remaining possibilities according to Table~\ref{TSubgroupsOfG2}.

Assume $H_1 = \SU(3)$. If $H_2 = \SU(3)$, then the $H'$-action on~$\LG_2$ is of
cohomogeneity one with principal isotropy group~$\Delta \SU(2)$,
see~\cite{hyperpolar}. If $H = \SU(3) \times \U(2)$, then an isotropy group of
the $H$-action on $\LG_2$ is $\Delta \U(2)$ whose slice representation is
equivalent to the standard representation of~$\SU(3)$ restricted to
$\SUxU{2}{1}$, which has trivial principal isotropy. If $H = \SU(3) \times
\LA_1^4$, then the principal isotropy is also trivial since a slice
representation is equivalent to the standard representation of~$\SO(3)$
on~$\C^3$. The action of $H = \SU(3) \times \SO(4)$ is of cohomogeneity
one~\cite{hyperpolar}. Let us show that the action of $H = \SU(3) \times
\widetilde{\U(2)}$ on~$\LG_2$ is of cohomogeneity two. The $7$-dimensional
irreducible representation of~$\LG_2$ splits as $\R^3 \oplus \R^4$ when
restricted to~$\widetilde{\U(2)}$, where $\widetilde{\U(2)}$ acts by the
standard $\U(2)$-representation on~$\R^4$ and by the adjoint representation
of~$\SU(2)$ on~$\R^3$. This representation, and hence the action of
$\widetilde{\U(2)}$ on~${\rm S}^6 = \LG_2 / \SU(3)$ has finite principal isotropy
groups. The action of $\SO(3)$ on ${\rm S}^6$ given by the $7$-dimensional
irreducible representation of~$\SO(3)$ has trivial principal isotropy groups
\cite{hh}, hence it follows that the principal isotropy groups of the $\SU(3)
\times \LA_1^{28}$-action on $\LG_2$ are trivial as well.

Now assume $H_1 = \LA_1^{28}$. Consider the action of $\LA_1^{28} \times
\LA_1^{28}$ on~$\LG_2$. Since by~\cite{wolfIrr} one slice representation is
equivalent to the $11$-dimensional irreducible representation of~$\SO(3)$, it
follows that the action has trivial principal isotropy groups~\cite{hh}.

Assume the action of $H' = \LA_1^{28} \times \SO(4)$ on~$\LG_2$ has principal
isotropy groups of positive dimension. Then it follows that the action is of
cohomogeneity at least six. We may assume that the identity element~$\e \in
\LG_2$ lies in a principal orbit. Then there is some non-zero element $X \in
\h' \cap \Delta \Lg_2$ acting trivially on the normal space $\N_{\e} \lp H'
\act \e \rp \subset \Lg_2$. But this contradicts the fact that the centralizer
of any non-zero element in~$\Lg_2$ is at most $4$-dimensional.
\qed\end{proof}


\section{Actions of cohomogeneity one or two}
\label{CohOneTwo}


In this section we classify all isometric actions of cohomogeneity less or
equal to two on the exceptional simple compact Lie groups $L = \LG_2$, $\LF_4$,
$\LE_6$, $\LE_7$ and $\LE_8$. By the results of~\cite{oniscik}, there are no
nontrivial transitive actions on the exceptional groups. In~\cite{hyperpolar},
a classification of cohomogeneity one actions on simple compact Lie groups was
obtained. However, these actions were only classified up to orbit equivalence
there and it remains to determine orbit equivalent subactions.

Obviously, we only need to consider such closed subgroups $H$ of $L \times L$
where $\dim H \ge \dim L - 2$. Since the lower bound on the dimension of groups
acting polarly on~$L$ given in Lemmata~\ref{LmDimSection} and \ref{LmDimBd} is
in all cases a weaker condition, all candidates for groups acting with
cohomogeneity one or two have already appeared in the proof of
Theorem~\ref{ThmPolar}. Indeed, in order to prove the classification theorems
~\ref{ThmCohOne} and \ref{ThmCohTwo} below, we will proceed in the same order
as in the proof of Theorem~\ref{ThmPolar}, i.e.\ we first consider subgroups of
$\Delta L$, then we study subactions of Hermann actions and finally we consider
all remaining subgroups of $L \times L$ which are of sufficient dimension.

The slice representations of an isometric Lie group action on a Riemannian
manifold of cohomogeneity one or two are orthogonal representations of
cohomogeneity one or two, respectively, and hence are polar. Therefore we may
immediately rule out all groups contained in a closed subgroup $H \subset L
\times L$ which has a non-polar slice representation. Here we can use the
classifications \cite{pth1}, \cite{hyperpolar} and \cite{polar} of
(hy\-per)\-po\-lar actions to exclude many candidates for cohomogeneity two
actions. In fact, the appearance of cohomogeneity two actions is one major
technical complication in~\cite{polar}, since for these actions the slice
representations do not contain any information about the polarity of the
action.


\begin{theorem}\label{ThmCohOne} Let $L$ be a simply connected
exceptional simple compact Lie group and let $H \subset L \times L$ be a closed
subgroup acting with cohomogeneity one on $L$. Then $H$ is conjugate to a group
$H_1 \times H_2$ or $H_2 \times H_1$ such that the triple $(L,H_1,H_2)$ occurs
in Table~\ref{TCohOne}. In particular, there are no isometric cohomogeneity one
actions on~$\LE_7$ and $\LE_8$ and any isometric cohomogeneity one action on
$\LE_6$ and $\LF_4$ is orbit equivalent to a Hermann action.
\end{theorem}

\begin{table}[!h]\rm
\begin{center}\begin{tabular}{|c|c|c|l|}
\hline \str $L$ & $H_1$ & $H_2$ & Description \\
\hline \hline

 $\LE_6$ & $\SU(6)\cdot\mySp(1)$ & $\LF_4$ & E\,II-IV \\

 $\LE_6$ & $\SU(6)\cdot\U(1)$ & $\LF_4$ &  \\

 $\LE_6$ & $\SU(6)$ & $\LF_4$ &  \\ \hline

 $\LE_6$ & $\Spin(10)\cdot\U(1)$ & $\LF_4$ & E\,III-IV \\ \hline \hline

 $\LF_4$ & $\mySp(3)\cdot\mySp(1)$ & $\Spin(9)$ & F\,I-II \\

 $\LF_4$ & $\mySp(3)\cdot\U(1)$ & $\Spin(9)$ &  \\

 $\LF_4$ & $\mySp(3)$ & $\Spin(9)$ &  \\ \hline

 $\LF_4$ & $\Spin(9)$ & $\Spin(9)$ & F\,II-II \\ \hline \hline

 $\LG_2$ & $\SU(3)$ & $\SU(3)$ &  \\ \hline

 $\LG_2$ & $\SU(3)$ & $\SO(4)$ &  \\ \hline

\end{tabular}
\caption{Cohomogeneity one actions.} \label{TCohOne}
\end{center}
\end{table}


\begin{theorem}\label{ThmCohTwo}
Let $L$ be a simply connected exceptional simple compact Lie group and let $H
\subset L \times L$ be a closed connected subgroup. Then the action of $H$ on
$L$ is of cohomogeneity two if and only if one of the following holds
\begin{enumerate}

\item $L = \LG_2$ and $H$ is conjugate to $\Delta \LG_2$;

\item $L = \LE_6$ and $H$ is conjugate to a group $$\lp \Spin(10) \times
    \Spin(10) \rp \cdot Q \subseteq \lp \Spin(10) \cdot \U(1) \rp \times
    \lp \Spin(10) \cdot \U(1) \rp,$$ where $Q \subset \U(1) \times \U(1)$
    is a one-dimensional closed connected subgroup such that $Q \neq \Delta
    \U(1)$.

\item $H$ is conjugate to a group $H_1 \times H_2$ or $H_2 \times H_1$ such
    that the triple $(L,H_1,H_2)$ occurs in Table~\ref{TCohTwo}.

\end{enumerate}
\end{theorem}

\paragraph{\bf Remarks.} For convenience, we have formulated the statement of
Theorems~\ref{ThmCohOne} and \ref{ThmCohTwo} only for simply connected groups
$L$; however, the result of course implies the classification also for the
non-simply connected case, since the cohomogeneity of an $H$-action on $L$
depends only on the conjugacy class of the subalgebra $\h \subset {\mathfrak l} \oplus
{\mathfrak l}$. In the last column of Tables~\ref{TCohOne} and \ref{TCohTwo}, the types of
the symmetric subgroups are given in case of a Hermann action, here we use the
notation as in Table~\ref{THermannExc}. Actions which appear in consecutive
rows of the tables without separating horizontal lines between them are orbit
equivalent to one another.

\begin{table}[!h]\rm
\begin{center}\begin{tabular}{|c|c|c|c|l|}
\hline \str $L$ & $H_1$ &  $H_2$ & polar? & Description \\
\hline \hline

 $\LE_6$ & $\mySp(4)/\{\pm1\}$ & $\Spin(10){\cdot}\U(1)$ & yes & E\,I-III \\ \hline

 $\LE_6$ & $\mySp(4)/\{\pm1\}$ & $\LF_4$ & yes & E\,I-IV \\ \hline

 $\LE_6$ & $\SU(6){\cdot}\mySp(1)$ & $\Spin(10){\cdot}\U(1)$ & yes & E\,II-III \\

 $\LE_6$ & $\SU(6){\cdot}\mySp(1)$ & $\Spin(10)$ &  &  \\ \hline

 $\LE_6$ & $\Spin(10){\cdot}\U(1)$ & $\Spin(10){\cdot}\U(1)$ & yes & E\,III-III \\ \hline

 $\LE_6$ & $\LF_4$ & $\LF_4$ & yes & E\,IV-IV \\ \hline

 $\LE_6$ & $\Spin(10)$ & $\LF_4$ & no &  \\ \hline

 $\LE_6$ & $\SUxU{5}{1}{\cdot}\mySp(1)$ & $\LF_4$ & no &  \\

 $\LE_6$ & $\SU(5){\cdot}\mySp(1)$ & $\LF_4$ & no &  \\ \hline \hline

 $\LE_7$ & $\SO'(12){\cdot}\mySp(1)$ & $\LE_6{\cdot}\U(1)$ & yes & E\,VI-VII \\

 $\LE_7$ & $\SO'(12){\cdot}\mySp(1)$ & $\LE_6$ &  &  \\ \hline \hline

 $\LF_4$ & $\Spin(9)$ & $\Spin(8)$ & no &  \\ \hline

 $\LF_4$ & $\Spin(9)$ & $\Spin(7)\cdot\SO(2)$ & no &  \\ \hline

 $\LF_4$ & $\Spin(9)$ & $\Spin(6)\cdot\Spin(3)$ & no &  \\ \hline

 $\LF_4$ & $\Spin(9)$ & $\LG_2^1\cdot\LA_1^8$ & no &  \\ \hline

 $\LF_4$ & $\Spin(9)$ & $\SU(3)\cdot\SU(3)$ & no &  \\ \hline \hline

 $\LG_2$ & $\SO(4)$ & $\SO(4)$ & yes & G \\ \hline

 $\LG_2$ & $\SU(3)$ & $\U(2)$ & no &  \\ \hline

 $\LG_2$ & $\SU(3)$ & $\widetilde{\U(2)}$ & no &  \\ \hline

\end{tabular}
\caption{Cohomogeneity two actions.} \label{TCohTwo}
\end{center}
\end{table}


\begin{proof}[Proof of Theorems~\ref{ThmCohOne} and \ref{ThmCohTwo}]
Let $H$ be a closed connected subgroup of a simply connected exceptional
compact Lie group~$L$ acting with cohomogeneity one or two. By
Lemma~\ref{LmMaxSubgr}, $H$ is either contained in $\Delta L$ or in $H_1 \times
H_2$, where $H_i \subset L$ are maximal connected subgroups. In the first case
it follows from the result of Section~\ref{Sigma} that the action is of
cohomogeneity~$\le 2$ if and only if $L = \LG_2$ and $H$ is conjugate
to~$\Delta \LG_2$.

\paragraph{\em Subactions of cohomogeneity one Hermann actions.}
Assume $H \subset H_1 \times H_2$, where $H_i \subset L$ are symmetric
subgroups such that the $H_1 \times H_2$-action on~$L$ is of cohomogeneity~$\le
1$. Since there are no transitive actions of this type~\cite{oniscik}, it
follows that the action is a Hermann action of cohomogeneity one, cf.\
Table~\ref{THermannExc}. (Note that the action F\,II-II does not appear in
Table~\ref{THermannExc}, since the groups $H_1$ and $H_2$ are conjugate.) We
will treat the subactions of these four actions in the following paragraphs.

\paragraph{\bf F\,II-II}
Assume there is a closed connected subgroup $H'$ of $H =
\Spin(9) \times \Spin(9)$ acting with cohomogeneity two on~$\LF_4$. The
$H$-orbit $H \act \e$ through the identity element~$\e \in \LF_4$ is the
subgroup $\Spin(9) \subset \LF_4$ and the slice representation of the isotropy
group $\Delta \Spin(9)$ at~$\e$ is equivalent to the $16$-dimensional spin
representation of~$\Spin(9)$. Consider the action of $H_{\e}' = H' \cap \Delta
\Spin(9)$ on the invariant subspace $\N_0 := \N_{\e} \lp H \act \e \rp$ of its
slice representation. Now there are two cases, depending on whether this action
is transitive on the sphere or not. If $H_{\e}'$ acts transitively on the unit
sphere in $\N_0$, then it follows that $\Delta \Spin(9)$ is contained in
$H_{\e}'$, since there is no non-trivial factorization~\cite{oniscik}
of~$\Spin(9)$. Since $\Delta \Spin(9) \subset \Spin(9) \times \Spin(9)$ is a
maximal connected subgroup, it follows that either $H'=\Delta \Spin(9)$, which
acts with cohomogeneity~$16$  or $H' = \Spin(9) \times \Spin(9)$. If $H_{\e}'$
does not act transitively on the unit sphere in $\N_0$, then it follows that
$H'$ acts transitively on $H \act \e$ and, again, since there is no non-trivial
factorization of $\Spin(9)$, it follows that $H'$ is of the form $\Spin(9)
\times K$ (or $K \times \Spin(9)$), where $K \subset \Spin(9)$ acts with
cohomogeneity two on the Cayley plane $\LF_4 / \Spin(9)$. These groups have
been classified in~\cite{pth1}, pp.~172-173.

\paragraph{\bf F\,II-I}
Assume $H'$ is a closed connected subgroup of $H = \Spin(9)
\times \lp \mySp(3) \cdot \mySp(1) \rp$ acting with cohomogeneity one or two
on~$\LF_4$. From the proof of Theorem~\ref{ThmSubCohOne} we see that it only
remains to consider the groups
\begin{align*}
\Spin(9) \times \lp \mySp(2) \cdot \mySp(1) \cdot \mySp(1) \rp,\quad
\Spin(9) \times \lp \mySp(3) \cdot \U(1) \rp,\quad
\Spin(9) \times \mySp(3).
\end{align*}
The first group acts with cohomogeneity greater than two~\cite{pth1},
pp.~172-173, the other groups are known to act with cohomogeneity
one~\cite{pth1}.

\paragraph{\bf E\,II-IV}
The groups considered in Theorem~\ref{ThmSubCohOne} cannot
act with cohomogeneity one or two since their dimension is too small. Thus it
remains to consider groups $H = H_1 \times H_2$ where either $H_1 = \SU(6)
\cdot \mySp(1)$ or $H_2 = \LF_4$.

Assume first $H_1 = \SU(6) \cdot \mySp(1)$. The argument in \cite{polar},
Section~12, p.~478, shows that there is no subgroup $H_2 \subsetneqq \LF_4$
acting with cohomogeneity two on $\LE_6$.

Assume now that $H_2 = \LF_4$. These actions have been studied in \cite{polar},
Section~12, pp.~477-478, it is shown there that some of them have non-polar
slice representations and thus are of cohomogeneity at least three;
furthermore, subactions of cohomogeneity one are determined. In the case where
$H_1 = \SUxU{5}{1} \cdot \mySp(1)$, it is (implicitly) shown in loc.\ cit.\ that
the action is of cohomogeneity two and that among its subactions $H_1 = \SU(5)
\cdot \mySp(1)$ is the only one of cohomogeneity two.  If $H = (\mySp(3) \cdot
\mySp(1)) \times \LF_4$ then we may assume $H_1 \subset H_2$ (see \cite{dynkin1},
p.~200) and the orbit through the identity element $\e \in \LE_6$ is the
subgroup~$\LF_4$. The slice representation at $\e$ is the $26$-dimensional
irreducible representation of~$\LF_4$ whose restriction to $\mySp(3) \cdot
\mySp(1)$ is non-polar by Proposition~\ref{PolMinCriteria}~(i) and
Table~\ref{TOrbitEqPolarSubGr}, hence not of cohomogeneity two. There are no
cohomogeneity two subactions since $\dim(H) = 76$.

\paragraph{\bf E\,III-IV}
Let $H$ be a closed connected subgroup of $(\Spin(10) \cdot
\U(1)) \times \LF_4$. Assume first $H$ contains the $\LF_4$-factor. Such
actions were studied in \cite{polar}, Section~12, pp.~478-479 where the
cohomogeneity two action of $\Spin(10) \times \LF_4$ is shown to be non-polar.
It is also shown in loc.\ cit.\ that all other actions of this type are
excluded by a dimension count or have non-polar slice representations. Assume
now $H$ contains the factor $\Spin(10) \cdot \U(1)$. By the argument in
\cite{polar}, Section~12, p.~479, we only have to consider the action of
$(\Spin(10) \cdot \U(1)) \times \Spin(9)$ which can easily be seen to have a
slice representation of cohomogeneity greater than two. In case $H$ does
contain neither the $\LF_4$-factor nor $\Spin(10) \cdot \U(1)$, it follows from
the last part of the proof of Theorem~\ref{ThmSubCohOne} that the action cannot
be of cohomogeneity two.

\paragraph{\em Subactions of cohomogeneity two Hermann actions.}
The orbit equivalent subactions of Hermann actions of cohomogeneity two were
determined in the proof of Theorem~\ref{ThmSubHermann}.

\paragraph{\em Other actions.}
It remains to determine those cohomogeneity two actions which are neither
subactions of Hermann actions nor of $\sigma$-actions.

\paragraph{\em Actions on $\LG_2$.}
For actions on $\LG_2$, the classification follows from
Theorem~\ref{ThmPrOrbG2}.

\paragraph{\em Actions on $\LF_4$.}
Adding up dimensions of the groups given in Table~\ref{TMaxF4}, we see that the
only actions of a group of sufficient dimension are
{\bf\ref{GrpF41}--\ref{GrpF43}} and {\bf\ref{GrpF41}--\ref{GrpF44}}. It is
shown in Section~\ref{PolarF4} that both actions are non-polar and of
cohomogeneity two and it follows from the arguments given there that there are
no proper subactions of cohomogeneity two.

\paragraph{\em Actions on $\LE_6$.}
Again by counting dimensions we only need to consider the action
{\bf\ref{GrpE64}--\ref{GrpE63}}. Since a slice representation of this action is
non-polar, it follows that the cohomogeneity is greater than two.

\paragraph{\em Actions on $\LE_7$.}
We only need to consider the action {\bf\ref{GrpE78}--\ref{GrpE75}}, which has
a non-polar slice representation.

\paragraph{\em Actions on $\LE_8$.}
Since the only closed connected subgroups in~$\LE_8$ of dimension~$\ge 110$ are
symmetric and Hermann actions on~$\LE_8$ are of cohomogeneity~$\ge 4$, we
conclude that there are no isometric actions on~$\LE_8$ of cohomogeneity one or
two.
\qed\end{proof}


\section{Low cohomogeneity actions on $\LE_8$}
\label{LowCoh}


\begin{theorem}\label{ThE8Low}
Let $H \subset \LE_8 \times \LE_8$ be a closed connected subgroup acting on
$\LE_8$. Then the $H$-action on $\LE_8$ is of cohomogeneity $k$ with $0 < k <
20$ if and only if it is conjugate to the action of one of the groups given in
Table~\ref{TE8Low}.
\end{theorem}

\begin{table}[!h]\rm
\begin{center}\begin{tabular}{|c|l||c|}
\hline \str  \makebox[115pt]{Subgroup of $\LE_8 \times \LE_8$} & Range & Cohomogeneity\ \\
\hline \hline

 $\Delta \LE_8$ &  & $8$  \\ \hline

 $\lp \LE_7 \times  \LE_7 \rp \cdot Q$ &$Q \subseteq \mySp(1) \times \mySp(1)$ & $10 - \dim Q$  \\ \hline

 $\SO'(16) \times \SO'(16)$ &  & $8$  \\ \hline

 $\lp \LE_7 \cdot P \rp \times \SO'(16)$ & $P \subseteq \mySp(1)$ & $\;\, 7 - \dim P$  \\ \hline

 $\lp \LE_7 \cdot P \rp \times \Spin(15)$ & $P \subseteq \mySp(1)$ & $10 - \dim P$  \\ \hline

\end{tabular}
\caption{Low cohomogeneity actions on $\LE_8$.} \label{TE8Low}
\end{center}
\end{table}


\begin{proof}
Assume $H \subset \LE_8 \times \LE_8$ is a closed connected subgroup acting
non-transitively and with cohomogeneity~$\le19$ on $\LE_8$. We know from
Lemma~\ref{LmMaxSubgr} that $H$ either is contained in $\Delta \LE_8$ or in a
group of the form $H \subseteq H_1 \times H_2$, where $H_i \subset \LE_8$ are
maximal closed connected subgroups. In the first case it follows that $H =
\Delta \LE_8$ since the maximal dimension of a proper closed subgroup
of~$\LE_8$ is~$136$. In the latter case, it follows that $\dim H_i \ge 93 =
\dim \LE_8 - \dim \LE_7 \cdot \mySp(1) -19 = 248 -136 -19$. We can see from
\cite{dynkin1} that the only such maximal connected subgroups of $\LE_8$ are
the two symmetric subgroups $\LE_7 \cdot \mySp(1)$ and $\SO'(16)$. The only
connected subgroups of these two groups whose dimension is greater than~$92$
are $\LE_7 \cdot \U(1)$, $\LE_7 \subset \LE_7 \cdot \mySp(1)$ and $\Spin(15)
\subset \SO'(16)$, respectively. Using the data of the slice representation of
the Hermann action E\,VIII-IX as given by Table~\ref{THermannSliceRep} and
looking at the isotropy representations of the symmetric spaces E\,VIII  and
E\,IX we can now verify the content of Table~\ref{TE8Low}.

Consider closed connected subgroups $H$ of $\lp \LE_7 \cdot \mySp(1) \rp \times
\lp \LE_7 \cdot \mySp(1) \rp$. Every such subgroup of dimension greater
than~$228$ is of the form $\lp \LE_7 \times \LE_7 \rp \cdot Q$, where $Q$ is a
closed connected subgroup of~$\mySp(1) \times \mySp(1)$. The principal isotropy
algebra for the action of $\lp \LE_7 \times  \LE_7 \rp \cdot Q$ on~$\LE_8$ is
isomorphic to $\spin(8)$ according to~\cite{hh}, p.~199; hence the
cohomogeneity of this action is $10 - \dim Q$. The cohomogeneity of the Hermann
action of $\SO'(16) \times \SO'(16)$ on~$\LE_8$ is equal to $\rk (\LE_8 /
\SO'(16)) = 8$. Now consider subactions of the Hermann action E\,VIII-IX. Every
subgroup of $\lp \LE_7 \cdot \mySp(1) \rp \times \SO'(16)$ of dimension greater
than~$228$ is of the form $\lp \LE_7 \cdot P \rp \times \SO'(16)$ or $\lp \LE_7
\cdot P \rp \times \Spin(15)$, where $P \subseteq \mySp(1)$ is a closed connected
subgroup. By Table~\ref{THermannSliceRep}, an isotropy subgroup of the Hermann
action E\,VIII-IX is locally isomorphic to $\U(8)$ and its slice representation
is on the Lie algebra level equivalent to the isotropy representation
of~$\SO(16) / \U(8)$. Following~\cite{hh}, the principal isotropy subalgebra is
$\su(2) \oplus \su(2) \oplus \su(2) \oplus \su(2)$. This shows that the
principal isotropy group is $12$-dimensional for all choices of~$P \subseteq
\mySp(1)$ since we have the inclusion $4 \cdot \su(2) \subset \su(8) \subset
\Le_7$. Finally, consider the action of $\lp \LE_7 \cdot \mySp(1) \rp \times
\Spin(15)$ on~$\LE_8$. It follows from the argument above that one isotropy
subgroup is locally isomorphic to~$\U(7)$. The isotropy representation of
$\SO(16) / \U(8)$ restricted to~$\U(7)$ splits into $\Lambda^2 \C^7 \oplus
\C^7$ and has finite principal isotropy subgroups. This shows that the actions
in the last row of Table~\ref{TE8Low} also have finite principal isotropy
subgroups.
\qed\end{proof}


\bibliographystyle{amsplain}

\begin{thebibliography}{XyzXyz2}

\bibitem{bergmann} I.\ Bergmann: {\it Reducible polar representations}, manuscripta
    math.\ {\bf 104}, 309-324 (2001)

\bibitem{biliotti} L.\ Biliotti: {\it Co\-iso\-tro\-pic and polar actions on compact
    irreducible Hermitian symmetric spaces}, Trans.\ Amer.\ Math.\ Soc.\ {\bf 358},
    3003-3022 (2006)

\bibitem{bg} L.\ Biliotti, A.\ Gori: {\it Co\-iso\-tro\-pic actions on complex
    Grassmannians}, Trans.\ Amer.\ Math.\ Soc.\ {\bf 357}, 1731-1751 (2005)

\bibitem{brueck} M.\ Br\"{u}ck: {\it \"{A}quifokale Familien in symmetrischen
    R\"{a}umen,
    (Equifocal families in symmetric spaces)}, (German) Doctoral Dissertation,
    Universit\"{a}t zu K\"{o}ln (1998)

\bibitem{dadok} J.\ Dadok: {\it Polar coordinates induced by actions of compact Lie
    groups}, Trans.\ Amer.\ Math.\ Soc.\ {\bf 288}, 125-137 (1985)

\bibitem{datri} J.E.\ D'Atri: {\it Certain isoparametric families of hypersurfaces
    in symmetric spaces}, J.\ Differential Geom.\ {\bf 14}, 21-40
    (1979)

\bibitem{dynkin1} E.B.\ Dynkin: {\it Semisimple subalgebras of the semisimple Lie
    algebras}, (Russian) Mat.\ Sbornik {\bf 30}, 349-462 (1952); English
    translation: Amer.\ Math.\ Soc.\ Transl.\ Ser.~2, {\bf 6}, 111-244 (1957)

\bibitem{dynkin2} E.B.\ Dynkin: {\it The maximal subgroups of the classical
    groups}, (Russian) Trudy Mosk.\ Mat.\ Obshch.\ {\bf 1} 39-166, (1952);
    English translation: Amer.\ Math.\ Soc.\ Transl.\ Ser.~2, {\bf 6} (1957)
    245-378, Zbl 0077.03403

\bibitem{eh2} J.-H.\ Eschenburg, E.\ Heintze: {\it On the classification of polar
    representations}, Math.\ Z.\ {\bf 232} (3) 391-398 (1999)

\bibitem{gorodski} C.\ Gorodski: {\it Polar actions on compact symmetric spaces
    which admit a totally geodesic principal orbit}, Geometriae Dedicata, {\bf
    103} (1), 193-204 (2004)

\bibitem{hl} E.\ Heintze, X.\ Liu: {\it A splitting theorem for isoparametric
    submanifolds in Hilbert space}, J.\ Differential Geom.\ {\bf 45}, 319-335
    (1997)

\bibitem{hptt} E.\ Heintze, R.\ Palais, C.-L.\ Terng, G.\ Thorbergsson:
    {\it Hyperpolar actions on symmetric spaces}, Geometry, topology and physics
    for Raoul Bott, ed.\ S.-T.\ Yau, International Press (1994)

\bibitem{helgason} S.\ Helgason: {\it Differential geometry, Lie groups and
    symmetric spaces}, Academic Press (1978)

\bibitem{hermann} R.\ Hermann: {\it Variational completeness for compact symmetric
    spaces}, Proc.\ Amer.\ Math.\ Soc.\ {\bf 11,} 544-546 (1960)

\bibitem{hh} W.C.\ Hsiang, W.Y.\ Hsiang: {\it Differentiable actions of compact
    connected classical groups II}, Ann.\ of Math.~(2) {\bf 92}, 189-223 (1970)

\bibitem{hsl} W.Y.\ Hsiang, H.B.\ Lawson: {\it Minimal submanifolds of low
    cohomogeneity}, J.\ Differential Geom.\ {\bf 5}, 1-38 (1971)

\bibitem{iwata} K.\ Iwata: {\it Compact transformation groups on rational
    cohomology Cayley projective planes}, T\^{o}hoku Math.\ J.\ {\bf 33},
    429-422 (1981)

\bibitem{hyperpolar} A.\ Kollross: {\it A classification of hyperpolar and
    cohomogeneity one actions}, Trans.\ Amer.\ Math.\ Soc.\ {\bf 354}, 571-612
    (2002)

\bibitem{polar} A.\ Kollross: {\it Polar actions on symmetric spaces}, J.\
    Differential Geom.\ {\bf 77} (3), 425-482 (2007)

\bibitem{kp} A.\ Kollross, F.\ Podest\`a: {\it Homogeneous spaces with polar
    isotropy}, ma\-nu\-scrip\-ta math.\ {\bf 110}~(4), 487-503 (2003)

\bibitem{oniscik} A.L.\ Oni\v{s}\v{c}ik: {\it Inclusion relations among transitive
    compact transformation groups}, (Russian) Trudy Mosk.\ Mat.\ Obshch.\ {\bf
    11}, 199-242, (1962); English translation: Amer.\ Math.\ Soc.\ Transl.\
    Ser.~2, {\bf 50}, 5-58 (1966)

\bibitem{oniscikBook} A.L.\ Oni\v{s}\v{c}ik: {\it Topology of transitive
    transformation groups}, Johann Ambrosius Barth (1994)

\bibitem{takagi} R.\ Takagi: {\it On homogeneous real hypersurfaces in a complex
    projective space}, Osaka J.\ Math.\ {\bf 10}, 495-506 (1973)

\bibitem{pt} R.\ Palais, C.-L.\ Terng: {\it A general theory of canonical forms},
    Trans.\ Amer.\ Math.\ Soc. {\bf 300}, 771-789 (1987)

\bibitem{pth1} F.\ Podest\`{a}, G.\ Thorbergsson: {\it Polar actions on rank one
    symmetric spaces}, J.\ Differential Geom.\ {\bf 53} (1), 131-175 (1999)

\bibitem{pth2} F.\ Podest\`{a}, G.\ Thorbergsson: {\it Polar and co\-iso\-tro\-pic
    actions on K\"{a}hler manifolds}, Trans.\ Amer.\ Math.\ Soc.\ {\bf 354},
    1749-1757 (2002)

\bibitem{wolfIrr} J.A.\ Wolf: {\it The geometry and structure of isotropy
    irreducible homogeneous spaces}, Acta Math.\ {\bf 120}, 59-148, (1968);
    correction, Acta Math.\ {\bf 152}, 141-142, (1984)

\end{thebibliography}

\end{document}